\documentclass[12pt,letterpaper]{article}
\usepackage{natbib}

\UseRawInputEncoding
\usepackage{bbm}
\usepackage{chngcntr}
\usepackage{multirow}

\usepackage[colorlinks=true,urlcolor=blue,citecolor=blue,linkcolor=blue,bookmarks=true,bookmarksopen=false]{hyperref}

\setcitestyle{round,aysep={,},yysep={;}, citesep={;}}

\usepackage{algorithmic}
\usepackage[linesnumbered,ruled,vlined]{algorithm2e}

\usepackage[left=1in,top=1in,right=1in,bottom=1in]{geometry}

\newcommand{\exclude}[1]{}

\renewcommand{\Omega}{\varOmega}
\newcommand{\R}{\mathbbm{R}}

\newcounter{commentcounter}
\long\def\symbolfootnote[#1]#2{\begingroup%
\def\thefootnote{\fnsymbol{footnote}}\footnote[#1]{#2}\endgroup}

\usepackage{sibarticle}

\usepackage{amsmath}
\usepackage{kbordermatrix}
\usepackage{mathtools}

\usepackage{tikz}
\usetikzlibrary{arrows,shapes,snakes,automata,backgrounds,petri}
\usetikzlibrary{graphs}

\newcommand{\ignore}[1]{}

\usepackage{color}
\usepackage{bm}
\usepackage{multirow}
\usepackage{booktabs}

\usepackage{rotating}

\usepackage{graphicx}
\usepackage{subfigure}
\usepackage{epstopdf}

\allowdisplaybreaks

\title{Mixed-Integer Programming for a Class of Robust Submodular Maximization Problems}
\ShortTitle{MIP for a Class of Robust Submodular Maximization Problems} \ShortAuthors{Authors}
\NumberOfAuthors{2} \FirstAuthor{Hsin-Yi Huang\footnotemark[1],  Hao-Hsiang Wu\footnotemark[1]  }
\FirstAuthorAddress{  Department of Management Science, National Yang Ming Chiao Tung University, Hsinchu, Taiwan\\ 
	{\tt	huanghsinyi.mg08@nycu.edu.tw, hhwu2@nycu.edu.tw }}

\SecondAuthor{Simge K\"u\c{c}\"ukyavuz\footnotemark[2]}
\SecondAuthorAddress{
	Department of Industrial Engineering and Management Sciences,
	Northwestern University, Evanston, IL {\tt
		simge@northwestern.edu}
}
\ShortAuthors{Huang, Wu, K\"u\c{c}\"ukyavuz}
\exclude{
	\SecondAuthor{}
	\SecondAuthorAddress{
		Department of Integrated Systems Engineering,
		The Ohio State University, {\tt
			wu.2294@osu.edu}
	}
}

\keywords{robust optimization; integer programming; facet; submodularity; sensor networks}

\begin{document}

\maketitle \centerline{\today}
\footnotetext[1]{Co-first authors ordered alphabetically.}
\footnotetext[2]{Corresponding author.} 

\begin{abstract}
We consider robust submodular maximization problems (RSMs), where given a set of $m$ monotone submodular objective functions, the robustness is with respect to the worst-case (scaled) objective function. The model we consider generalizes two variants of robust submodular maximization problems in the literature, depending on the choice of the scaling vector. On one hand, by using unit scaling, we obtain a usual robust submodular maximization problem. On the other hand, by letting the scaling vector be the optimal objective function of each individual (NP-hard) submodular maximization problem, we obtain a second variant. While the robust version of the objective is no longer submodular, we reformulate the problem by exploiting the submodularity of each function. We conduct a polyhedral study of the resulting formulation and provide conditions under which the submodular inequalities are facet-defining for a key mixed-integer set.  We investigate several strategies for incorporating these inequalities within a delayed cut generation framework to solve the problem exactly. For the second variant, we provide an algorithm to obtain a feasible solution along with its optimality gap. We apply the proposed methods to a sensor placement optimization problem in water distribution networks using real-world datasets to demonstrate the effectiveness of the methods.
\end{abstract}

\renewcommand{\baselinestretch}{1.5}

\section{Introduction}\label{problem}

We study two variants of robust submodular  maximization problems (RSMs) considered in \cite{RSM2008} and \cite{RIM2016}, where the robustness is with respect to the worst case of a finite number of (scaled) submodular functions.  Specifically, let $V = \{1,\dots,n\}$ be a finite non-empty ground set, where $n \in \mathbb {N}$. Let $[m]=\{1,\dots,m\}$ be the set of the first $m \in \mathbb{N}$ positive integers.  For all $i \in [m]$, a function $f_i:2^V \rightarrow \mathbb R$ is submodular if 
\begin{equation*}
 f_i(X \cup \{j\})-f_i(X) \ge f_i(X' \cup \{j\})-f_i(X') \text{ for } X' \subseteq X \subseteq V \text{ and } j \in V \setminus X.
\end{equation*} 
This   definition of submodularity uses the concept of a  marginal contribution. In particular, the term $f_i(X \cup \{j\})-f_i(X)$ denotes the marginal contribution of the element $j$ when added to the set $X$ in function $f_i$, and the marginal contribution of $j$ decreases if the set $X$ includes more elements from the set $V\setminus X$. Given    monotonically non-decreasing submodular functions, $f_i$, we assume, without loss of generality, that $f_i(\emptyset) = 0, i \in [m]$.  Note that, throughout the paper, we  use the notation $\bar {\bold{x}} \in \mathbb B^{n}$ and its support $\bar X  = \{i \in V: \bar x_i = 1 \}$, and refer to the corresponding function evaluations $f_i(\bar {\bold{x}})$ for  $\bar {\bold{x}} \in \mathbb B^{n}$ and $f_i(\bar X)$ for the corresponding support $\bar X \subseteq {V}$,  interchangeably. 
Let $\mathcal{X}$ be a set of constraints on the binary variables $\bold{x}  \in \mathbb B^{n}$. 

Given a single monotone submodular set function $f_i(\cdot)$, the traditional  submodular maximization problem is defined as
\begin{equation}\label{SubMax}
	\max_{\bold{x} \in \mathcal{X}\cap \mathbb B^{n}}f_i(\bold{x}).
\end{equation} 
It is well-known that submodular maximization is NP-hard. 

 \cite{RSM2008} study a robust variant  of Problem \eqref{SubMax}, where given $m$ submodular functions $f_i: 2^V \rightarrow \mathbb R, i\in[m]$, the objective is to maximize  the worst case (minimum) of these $m$ submodular functions, i.e., 
\begin{equation}\label{RSM2}
	\max_{\bold{x} \in \mathcal{X}}  \min_{i \in [m]} f_i(\bold{x}).
\end{equation} 
In other words, Problem \eqref{RSM2} aims to find a solution $\bold{x} \in \mathcal X$ that is robust against the minimum possible value given by $\min_{i \in [m]} f_i(\bold{x})$. That is, an optimal solution $\bold{x^*} \in \mathcal X \cap \mathbb B^{n}$ satisfies $\min_{i \in [m]} f_i(\bold{x^*}) \ge \min_{i \in [m]} f_i(\bold{\bar x})$ for all $\bold{\bar x} \in \mathcal X \cap \mathbb B^{n}$. Problem \eqref{RSM2}, introduced by \cite{RSM2008},  is the  first robust extension of submodular maximization,  and it inspired various extensions of  robustness such as \cite{RIM2016,Bogunovic2017,Orlin2018,Staib2019,Adibi2022}.

In this paper, in addition to the basic RSM Problem \eqref{RSM2}, we also consider  the formulation of \cite{RIM2016}, which extends the robustness of Problem \eqref{RSM2} to consider the performance of the robust solution in proportion to the performance of the optimal solution for each submodular function. More precisely,  let $\bold{x}^*_i$ be an optimal solution of the $i$-th traditional submodular maximization problem \eqref{SubMax}. The RSM of \cite{RIM2016} is defined as
\begin{equation}\label{RSM3}
	\max_{\bold{x} \in \mathcal{X}\cap \mathbb B^{n}}  \min_{i \in [m]} \frac{f_i(\bold{x})}{f_i(\bold{x}_i^*)}.
\end{equation} 
For $\bold{ x} \in \mathcal X$, the authors consider the proportion of the  function value $f_i(\bold{x})$ to the largest possible function value $f_i(\bold{x}_i^*)$ for each $i \in [m]$. Problem \eqref{RSM3} aims to find a solution $\bold{ x} \in \mathcal X$ that maximizes the worst (smallest) value of these $m$ proportions. In other words, the optimal solution $\bold{x^*}$ of Problem \eqref{RSM3} satisfies $\min_{i \in [m]} \frac{f_i(\bold{x^*})}{f_i(\bold{x}_i^*)} \ge \min_{i \in [m]} \frac{f_i(\bold{\bar x})}{f_i(\bold{x}_i^*)}$ for all $\bold{\bar x} \in \mathcal{X}$.

In fact, we observe that Problems \eqref{RSM2} and \eqref{RSM3} can be generalized as the problem 
\begin{equation}\label{RSM}
\max_{\bold{x} \in \mathcal{X}\cap \mathbb B^{n}}  \min_{i \in [m]} \frac{f_i(\bold{x})}{\alpha_i},
\end{equation}
where  $\boldsymbol{\alpha}=(\alpha_1, \alpha_2, \dots, \alpha_m) \in \mathbb R^m_+$ is a given vector of nonnegative scalars. 
Problem  \eqref{RSM} is equivalent to Problem \eqref{RSM2} under the case $\boldsymbol{\alpha} = \bold{1}$. Furthermore, if we  solve $m$ submodular maximization problems and let $\boldsymbol{\alpha} = (f_1(\bold{x}_1^*),f_2(\bold{x}_2^*), \dots, f_m(\bold{x}_m^*))$, then Problem  \eqref{RSM3} is the same as Problem  \eqref{RSM}.

\cite{RSM2008} review a wide range of  applications of RSMs. For example, sensor placement optimization  for detecting the contamination of water networks \cite[]{KrauseWater2008,Leskovec07,RSM2008} can be modeled in the form of Problem \eqref{RSM}. Note that for this application in critical infrastructure we must take into account the issues of public health and security \cite[]{BWSN2008}, and rather than a placement that optimizes an expected performance measure, we are interested in optimizing the worst-case performance.  
Under public health considerations, a relevant objective concerns the population affected by the pollutant, where either the exact number or the proportion of people far from the pollutant is relevant. For example, functions $f_1(\bold{\bar x}) = 2$ and $f_2(\bold{\bar x}) = 10$ capture the exact number of individuals protected from the outbreak by decision $\bar x$ under $m = 2$ scenarios. Using the first performance measure, $f_1(\bold{\bar x}) = 2$ is the worst case of two scenarios. However, if a decision maker initially assesses an ability to protect $\alpha_1 = 2$ and $\alpha_2 = 20$ individuals for the first and second scenarios respectively, then the second scenario with $\frac{f_2(\bold{\bar x})}{\alpha_2} = 0.5$ has the worst proportion compared to $\frac{f_1(\bold{\bar x})}{\alpha_1} = 1$. A higher value of $\alpha_i$ for all $i \in [m]$ indicates an ambition to protect more individuals in a given scenario; however, because of the limitation of resources, the largest number of saved individuals  cannot be greater than $f_i(\bold{x}_i^*)$ for all scenarios $i \in [m]$. Therefore, it is reasonable to assume that $1 \le \alpha_i \le f_i(\bold{x}_i^*)$ for all $i \in [m]$. In our computational study, we demonstrate the effectiveness of our proposed methods on this sensor placement optimization problem. The detailed model of \cite{KrauseWater2008,Leskovec07} for Problem \eqref{RSM} will be given in Section \ref{computation}.

Previous literature on RSM focuses on a bicriteria approximation of the relaxation of Problem \eqref{RSM2} under certain constraints. \cite{RSM2008} show that under the cardinality constraint $\mathcal{X}_c= \{x: \sum_{i \in V } x_i \le b\}$ and $b \in \mathbb N$, there is no constant-ratio approximation algorithm for solving Problem \eqref{RSM2} unless $\mathcal{NP}=\mathcal{P}$. \cite{RSM2008} propose the SATURATE algorithm that provides a solution $\bold{\bar x}_s$ such that the objective value $\min_{i \in [m]} f_i(\bold{\bar x}_s) \ge \max_{\bold{x} \in \mathcal{X}_c \cap \mathbb B^{n}}  \min_{i \in [m]} f_i(\bold{x})$, where ${||\bold{\bar x}_s||}_0 \le \lambda_s b$ and $\lambda_s = 1+ \log (\max_{j\in V} \sum_{i \in [m]}f_i(\{j\}))$. \cite{Powers2016} subsequently propose the GENSAT algorithm {under an assumption that the submodular maximization problem with a matroid constraint has an approximation guarantee, $\lambda_g$.}
 For a fixed $\tau \in \mathbb R$, given a $\beta \in \mathbb R$, GENSAT provides a lower bound $\beta \tau$ for the minimal value of  every fraction $\gamma$ of $m$ submodular functions, where $\gamma \ge \frac{\lambda_g-\beta}{1-\beta}$ 
{and $\lambda_g \in \mathbb R$ is an  approximation guarantee based on the assumption shown in  Theorem 1 of \cite{Powers2016}}.  For the case  $\boldsymbol{\alpha} = (f_1(\bold{x}_1^*),f_2(\bold{x}_2^*), \dots, f_m(\bold{x}_m^*))$ in Problem \eqref{RSM} under cardinality constraint $\mathcal{X}_c$, \cite{RIM2016} show a strong approximation hardness result that the bicriteria approximation has to select at least a factor of $\lceil b \log m \rceil$ elements from $V$. Despite the hardness of solving the RSMs shown in \cite{RSM2008,RIM2016}, our research interest is to study the mathematical structure of Problem \eqref{RSM}. Instead of the approximation methods, the main goal of this paper is to provide exact methods based on mixed-integer programming  and polyhedral theory to solve Problem \eqref{RSM}, leveraging the tremendous power of mixed-integer programming solvers in obtaining solutions to many NP-hard problems.

Numerous optimization problems involving submodularity have been investigated via a mixed--integer programming lens, including but not limited to,  submodular  maximization  \cite[]{NW81,Ahmed2011,first2016,Yu2017,Shi2022,Coniglio2022},  submodular minimization  \cite[]{YuQ2022,QYu2023},  conic  quadratic optimization  \cite[]{Andres2018,Alper2020,Alper2022,kilincc2020conic},  $k$-submodular optimization  \cite[]{QYuk2021,QYu2021}, and  chance-constrained optimization  \cite[]{Wu2017, Fatma2022, Shen2022}.  We refer the reader to a recent tutorial \cite[]{KYTut2023} for an overview of these approaches.  Motivated by the success of these approaches in finding exact solutions to challenging submodular optimization problems, in this paper, we also undertake a polyhedral approach for Problem \eqref{RSM}, which is a robust version of the submodular maximization problem \eqref{SubMax}. One immediate challenge we face, as we will see later, is that the robust objective is no longer submodular even if each individual function is submodular.

{Robust optimization aims to deal with the worst-case over uncertain data with a broad array of applications such as finance \cite[]{Ghaoui2003,Goldfarb2003,Tutuncu2004}, supply chain management \cite[]{BenTalSupply2005,BertsimasSupply2006}, social networks \cite[]{RIM2016,Nannicini2019}, and energy systems \cite[]{Mulvey1995,Zhao2012,Bertsimas2013}. We refer the reader to the survey of \cite{Kouvelis1997,Bertsimas2011} for an overview of various domains. There are scalable algorithms for robust convex optimization  \cite[]{Ben1998,Ben1999,Ben2000}, robust discrete optimization under certain uncertainty sets \cite[]{Bertsimas2003,Bertsimas2004,Alper2006}, and two-stage robust linear programming \cite[]{Zhao2012,Jiang2012,Bertsimas2013,Zeng2013}, mainly relying on duality results of convex (or linear) programs. 
 However, submodular functions are neither convex nor concave, in general. Therefore these approaches are not directly applicable for the robust submodular optimization problem we consider.

Recall that Problem \eqref{RSM} is a robust version of the submodular maximization problem \eqref{SubMax}. Given $i \in [m]$, Problem \eqref{SubMax} is a class of $\mathcal {NP}$-hard problems  \cite[see, e.g.,][]{Feige1998, Feige2011}. In addition to network optimization \cite[]{MAXCOV74,KKT03,first2016,FischettiSocial2018,Cordeau2019,Dilek2019}, submodular maximization appears in other  modern applications including but  not limited to public security and health \cite[]{Leskovec07, KrauseWater2008, ZhengIISE2019}, computer vision \cite[]{Boykov2001,Jegelka2011}, computational linguistics \cite[]{Lin2011}, and artificial intelligence \cite[]{KrauseSensor2008,Golovin2011}. We  refer the reader to the survey of \cite{SubMaxApp2012} for an overview of various application domains of submodular optimization. There are two well-known approaches for solving Problem \eqref{SubMax}, either  exactly using delayed constraint generation approaches  or approximately using the greedy method  based on the seminal results of \cite{NW81} and  \cite{NWF78}, respectively. The greedy method has  $(1-1/e)$ optimality guarantee for monotone submodular maximization under a cardinality constraint  $\mathcal{X}_c$. 
 For a stochastic (expected value) version of Problem \eqref{SubMax} with a finite number of scenarios, \cite{first2016} introduce a two-stage stochastic submodular optimization model assuming that the second-stage objective function is submodular, where a corresponding delayed constraint generation algorithm with the submodular inequality of \cite{NW81} can be used for solving the problem. The expectation of stochastic submodular functions preserves submodularity, thereby enabling the adaptation of methods that exploit submodularity to the stochastic case.

 In contrast, in this paper, we consider a robust variant of monotone submodular function maximization (Problem \eqref{RSM}).  There are three difficulties with solving Problem \eqref{RSM}. First, for a given $\bold{x} \in \mathcal{X}$, the objective $\min_{i \in [m]} \frac{f_i(\bold{ x})}{\alpha_i}$ loses the submodularity property, and one cannot use the method of \cite{NW81}  directly. Second, we do not restrict ourselves to a particular type of constraint set (such as cardinality) in  $\mathcal{X}$, therefore any algorithm that assumes a particular constraint structure cannot be immediately applied. Finally, under the special case $\boldsymbol{\alpha} = (f_1(\bold{x}_1^*),f_2(\bold{x}_2^*), \dots, f_m(\bold{x}_m^*))$, it is very hard to solve  $m$ $\mathcal{NP}$-hard problems within a reasonable period of an execution time limit in order to define Problem \eqref{RSM3}. To conquer these difficulties, we provide an alternative formulation of Problem \eqref{RSM} that allows us to leverage the known submodular inequalities. We then conduct a polyhedral study of the associated mixed-integer set. Finally, for the hard special case with $\boldsymbol{\alpha} = (f_1(\bold{x}_1^*),f_2(\bold{x}_2^*), \dots, f_m(\bold{x}_m^*))$, we provide an algorithm that obtains a near-optimal solution equipped with an optimality gap.

The contributions and the outline of this paper are summarized as follows. In Section \ref{section_methods}, we review an alternative piecewise-linear  reformulation of Problem \eqref{RSM}, which enables the use of the submodular inequalities of \cite{NW81}. We conduct a polyhedral analysis of the associated mixed-integer set given by the alternative formulation and propose a facet-defining condition for the submodular inequalities. For the special case of Problem \eqref{RSM3}, we propose a method to estimate the optimality gap of the problem if it is too time-consuming to  obtain the optimal value of $\alpha_{i} = f_i(\bold{x}_i^*)$ for all $i \in [m]$. Based on these analyses, we  investigate several computational strategies   and propose a delayed constraint generation algorithm for Problem \eqref{RSM}. Finally, in Section \ref{computation}, we demonstrate the proposed methods on a sensor placement optimization problem in water networks using real-world datasets. We conclude in Section \ref{sec:conc}.

\section{Models and Methods} \label{section_methods}
In this section, we investigate models and methods for  Problem \eqref{RSM}. \cite{RSM2008} observe that the objective $\min_{i \in [m]} f_i(\bold{ x})$ of Problem \eqref{RSM2} is no longer submodular, even though each individual function $f_i$ is submodular. Therefore, Problem \eqref{RSM} also loses the submodularity property in the associated objective $\min_{i \in [m]} \frac{f_i(\bold{ x})}{\alpha_i}$ even for the case $\boldsymbol{\alpha} = \bold{1}$. However, we propose  an alternative formulation that exploits the submodularity property of each individual function. This alternative formulation is  crucial to derive several approaches to solve Problem \eqref{RSM}.

\subsection{An Alternative Formulation}

We first consider the alternative formulation of Problem \eqref{RSM}. Given constants $\alpha_i, i\in[m]$, the formulation is defined as
\begin{subequations}\label{RSM_alt}
	\begin{align}
	\max~~& \eta \\
	\text{s.t.}~~
	& \eta \le \frac{\theta_i}{\alpha_i} & \forall i \in [m] \label{RSM_alt:const}\\      
	& \theta_i \le f_i(\bold{x})  & \forall i \in [m] \label{RSM_alt:const2}\\      
	&\bold{x} \in \mathcal{X}\cap \mathbb B^{n}, \eta \in \mathbb R, \boldsymbol{\theta} \in \mathbb R^{m},\label{RSM_alt:const3}
	\end{align} 
\end{subequations}
where $\eta \in \mathbb R$ is a variable that captures the value of $\min_{i \in [m]} \frac{f_i(\bold{x})}{\alpha_i}$, and $\boldsymbol{\theta}$ is an $m$-dimensional vector of variables $\theta_i$ lower bounding the value of $f_i(\bold{x})$ for each $i \in [m]$. 
Note that in Formulation \eqref{RSM_alt},  constraints \eqref{RSM_alt:const2} entail the hypograph of $m$ submodular functions.  Since the function $f_i(\bold{x})$ of Formulation \eqref{RSM_alt} is submodular over the domain $\mathbb B^n$ for all $i \in [m]$, its hypograph is defined by  submodular inequalities of \cite{NW81}, given by  
\begin{equation}\label{sub_cut}
\theta_i \le f_i(S) - \sum_{j\in S}  \rho^i_j(V\setminus\{j\})(1-x_j)+\sum_{j\in V\setminus S} \rho^i_j(S) x_j , \forall S\subseteq V,
\end{equation} 
where $\rho_j^i(S)=f_i(S\cup\{j\})-f_i(S)$ captures the marginal contribution of including $j\in V\setminus S$ to a subset $S$. 
Using this observation, we derive a mixed-integer linear programming reformulation, where constraint \eqref{RSM_alt:const2} is replaced by inequalities \eqref{sub_cut} for all $i\in[m]$. Furthermore, the variables $\theta_i, i\in[m]$ can be projected out to arrive at the formulation
\begin{subequations}\label{RSM_alt2}
	\begin{align}
	\max~~& \eta \\
	\text{s.t.}~~
	& \eta \le 
	\frac{1}{\alpha_i}( f_i(S) - \sum_{j\in S}  \rho^i_j(V\setminus\{j\})(1-x_j)+\sum_{j\in V\setminus S} \rho^i_j(S) x_j) , \forall S\subseteq V, i\in[m] \label{eq:submod-eta}\\
	&\bold{x} \in \mathcal{X}\cap \mathbb B^{n}, \eta \in \mathbb R.
	\end{align} 
\end{subequations}	
The resulting formulation \eqref{RSM_alt2} has exponentially many constraints. Hence, we 
propose a delayed constraint generation (DCG) method to solve Formulation \eqref{RSM_alt}. In the proposed model, a relaxed master problem (RMP) at any iteration is formulated as
\begin{subequations}\label{RMP}
	\begin{align}
	\max~~& \eta \\
	\text{s.t.}~~
	& (\eta, \bold{x}) \in \mathcal{C} \label{RMP2}\\      
	&\bold{x} \in \mathcal{X}\cap \mathbb B^{n}, \eta \in \mathbb R,
	\end{align} 
\end{subequations} 
where {$\mathcal{C}$ is a mixed-integer set defined by a subset of the constraints \eqref{eq:submod-eta}} generated until the current iteration. In the next subsection, we  consider how to choose the inequalities to include in the set $\mathcal C$.

\subsection{Analyses of the Submodular Inequality for  RSM}

First, we observe that for large $m$, adding a submodular inequality \eqref{eq:submod-eta} for each  $i$ in a DCG algorithm may be inefficient. Motivated by this, we make a key observation that  a  mixed-integer set that includes fewer submodular inequalities compared to the submodular inequalities for all $i \in  [m]$ is sufficient to define $\mathcal{C}$ to find an optimal solution of Problem \eqref{RSM}. Before we give our analysis, we provide some useful definitions that identify an important index that determines the minimum of $m$ submodular functions for a given set.

\begin{definition}
	Given a subset $S \subseteq V$, we define a function 
\begin{equation*}
	\bold i(S) = \mathop{\arg\min}_{i \in [m]} \frac{f_i(S)}{\alpha_i},
\end{equation*}
where the function $\bold i:2^V \rightarrow \mathbb N$ returns the value of $i$ for which $\frac{f_i(S)}{\alpha_i}$ is the smallest. In other words, given a subset $S \subseteq V$, the corresponding value $\bold i(S)$ denotes an index such that $\frac{f_{\bold i(S)}(S)}{\alpha_{\bold i(S)}} \le \frac{f_i(S)}{\alpha_i}$ for all $i \in [m]$.
\end{definition}
Throughout this paper, the function $\bold i$ plays a key role in providing an upper bound for Problem \eqref{RSM}. Based on this index function, we define a mixed-integer set $\mathcal F$ as  
\begin{equation}\label{MIP_F}
\mathcal F=\{(\eta,\bold{x})\in \R \times \mathbb B^{n}:\eta \le \frac{1}{\alpha_{\bold i(S)}} ( f_{\bold i(S)}(S) - \sum_{j\in S}  \rho^{\bold i(S)}_j(V\setminus\{j\})(1-x_j)+\sum_{j\in V\setminus S} \rho^{\bold i(S)}_j(S) x_j), \forall S\subseteq V\}.
\end{equation}

In what follows, we prove that we can let $\mathcal{C}=\mathcal F$ in constraint \eqref{RMP2} of  RMP \eqref{RMP}.
\begin{proposition}\label{prop:MIP_valid}
The mixed-integer set $\mathcal F$ is sufficient for defining  $\mathcal{C}$ in  \eqref{RMP2} of  RMP \eqref{RMP}  to find an optimal solution of Problem \eqref{RSM}.
\end{proposition}    
\begin{proof}
\cite{NW81} show the validity of submodular inequality \eqref{sub_cut}. Thus, we have
	\begin{align*}
	\eta & \le \frac{\theta_{\bold i(S )}}{\alpha_{\bold i(S )}}  \le   \frac{f_{\bold i(S )}(S )}{\alpha_{\bold i(S )}}\\
	& \le \frac{ f_{\bold i(S )}(S ) - \sum_{j\in S }  \rho^{\bold i(S )}_j(V\setminus\{j\})(1-x_j)+\sum_{j\in V\setminus S } \rho^{\bold i(S )}_j(S ) x_j}{\alpha_{\bold i(S )}}.
	\end{align*} 
Therefore, Problem \eqref{RSM} is equivalent to the mixed-integer linear program
	\begin{align*}
	\max~~& \eta \\
	\text{s.t.}~~
	& \eta \le \frac{\theta_i}{\alpha_i} & \forall i \in [m] \\      
	& \theta_{\bold i(S )} \le f_{\bold i(S )}(S ) - \sum_{j\in S }  \rho^{\bold i(S )}_j(V\setminus\{j\})(1-x_j)+\sum_{j\in V\setminus S } \rho^{\bold i(S )}_j(S ) x_j  & \forall S \subseteq V\ \\      
	&\bold{S } \in \mathcal{X}\cap \mathbb B^{n}, \eta \in \mathbb R, \boldsymbol{\theta} \in \mathbb R^{m}.
	\end{align*} 
Projecting out the $\theta$ variables, we obtain the desired result.
\end{proof} 

Proposition \ref{prop:MIP_valid} shows that given $S  \subset V$, it is sufficient to add a submodular inequality \eqref{eq:submod-eta} with  $i = \bold i(S )$ to define the set $\mathcal C$ in  RMP \eqref{RSM} for solving the problem. Note that considering  all submodular inequalities defining $F_i, i \in  [m]$ may give a stronger formulation than considering $\mathcal F$. In our computational experiments, we observe that obtaining a violated submodular inequality is time-consuming, and as such, Proposition \ref{prop:MIP_valid} plays an important role in reducing the total number of inequalities. Below, we provide further analysis of $\mathcal F$ to improve algorithmic efficiency.

We start by providing a proposition that gives sufficient conditions under which the submodular inequality \eqref{eq:submod-eta} is facet-defining for conv($\mathcal F$). Let $\bold {e}_j$ be the $j$th unit vector of appropriate dimension. 
\begin{proposition}\label{prop:sub_facet}
Given  $S  \subseteq V$ and $\bar i \in [m]$, the submodular inequality
\begin{equation*}
	\eta \le \frac{1}{\alpha_{\bar i}} (f_{\bar i}(S ) - \sum_{j\in S }  \rho^{\bar i}_j(V\setminus\{j\})(1-x_j)+\sum_{j\in V\setminus S } \rho^{\bar i}_j(S ) x_j)
\end{equation*} 
is facet defining for conv($\mathcal F$) if the following  conditions hold:
	\begin{itemize}		
	\item[(i)] for any $j \in S $, there  exists at least one element $k_j \in V \setminus S  $ such that $\rho^{\bar i}_j(\{k_j\}) = 0$ and   $\frac{f_{\bar i}(S )}{\alpha_{\bar i}} = \frac{f_{\bold i(S )}(S )}{\alpha_{\bold i(S )}}
		= \frac{f_{\bold i(S  \setminus \{j\} \cup \{k_j\})}(S  \setminus \{j\} \cup \{k_j\})}{\alpha_{\bold i(S  \setminus \{j\} \cup \{k_j\})}}
		= \frac{f_{\bold i(S  \cup \{k_j\})}(S  \cup \{k_j\})}{\alpha_{\bold i(S  \cup \{k_j\})}}$	
	\item[(ii)]  {for any $j \in V\setminus S $, we have $ \frac{f_{\bar i}(S )+\rho^{\bar i}_j(S )}{\alpha_{\bar i} }= \frac{f_{\bold i(S  \cup \{j\})}(S  \cup \{j\})}{\alpha_{\bold i(S  \cup \{j\})}} $, where $\bar i = \bold i(S )$.} 
\end{itemize}
\end{proposition} 
\begin{proof}
Since dim($\mathcal F$) = $n+1$, we enumerate $n+1$ affinely independent points on the face defined by the submodular inequality \eqref{eq:submod-eta} under conditions (i) and (ii). 

	\begin{itemize}		
	\item[(a)] Given  $S  \subseteq V$, consider the point ($\eta, \bold{x}$) = {($\frac{f_{\bold i(S )}(S )}{\alpha_{\bold i(S )}},  \sum_{i \in S } \bold {e}_i$)} on the face defined by  inequality \eqref{eq:submod-eta}. 
	\item[(b)]  Building on the point given in  (a), we consider a set of points $P$, where $|P| = | V \setminus S |$ and each point $(\eta, \bold{x}) \in P$ is given by {$(
	\frac{f_{\bold i(\bold{S  \cup \{j\}})}(S  \cup \{j\})}{\alpha_{\bold i(S  \cup \{j\})}}, \sum_{i \in S } \bold {e}_i + \bold{e}_j)$} for all $j \in V\setminus S $, which is on the face defined by  inequality \eqref{eq:submod-eta} under condition (ii).  
	\item[(c)] From conditions (i) and (ii), for any $j \in S $, there exists  $k_j \in V \setminus S  $ such that $\rho^{\bar i}_j(\{k_j\}) = 0$ and $\bar i = \bold i(S )$. We conclude that $\rho^{\bold i(S )}_j(V\setminus\{j\}) = 0$ for any $j \in S $.   {Note that $ \frac{f_{\bold i(S  \setminus\{j\} \cup \{k_{j}\})}(S  \setminus\{j\} \cup \{k_{j}\})}{\alpha_{\bold i(S  \setminus\{j\} \cup \{k_{j}\})}} = \frac{f_{\bold i(S  \cup \{k_j\})}(S   \cup \{k_{j}\})-\rho^{\bold i(S  \cup \{k_j\})}_j(S  \cup \{k_j\})}{\alpha_{\bold i(S  \cup \{k_j\})}} = \frac{f_{\bold i(S  \cup \{k_j\})}(S   \cup \{k_{j}\}) + 0}{\alpha_{\bold i(S  \cup \{k_j\})}}$ since $\rho^{\bold i(S  \cup \{k_j\})}_j(S  \cup \{k_j\}) = \rho^{\bar i}_j( \{k_j\}) = 0$ and $\frac{f_{\bar i}(S )}{\alpha_{\bar i}} = \frac{f_{\bold i(S )}(S )}{\alpha_{\bold i(S )}}
		= \frac{f_{\bold i(S  \setminus \{j\} \cup \{k_j\})}(S  \setminus \{j\} \cup \{k_j\})}{\alpha_{\bold i(S  \setminus \{j\} \cup \{k_j\})}}
		= \frac{f_{\bold i(S  \cup \{k_j\})}(S  \cup \{k_j\})}{\alpha_{\bold i(S  \cup \{k_j\})}} $. Therefore, we obtain a  set of points $\bar P$ on the face defined by  inequality \eqref{eq:submod-eta}, where $(\eta, \bold{x})=(\frac{f_{\bold i(S  \setminus\{j\} \cup \{k_{j}\})}(S  \setminus\{j\} \cup \{k_{j}\})}{\alpha_{\bold i(S  \setminus\{j\} \cup \{k_{j}\})}}, \sum_{i \in S  \setminus\{j\}} \bold {e}_i + \bold{e}_{k_{j}}) \in \bar P$ for all $j \in S $ and $|\bar P| = |S |$.	  }

\end{itemize}	
	Note that these $n+1$ points can be represented as an $(n+1)\times(n+1)$ matrix, where the first $| V \setminus S |$ rows are the points  $P$ described in (b), from the $(| V \setminus S |+1)$-th row to the $|V|$-th row are the points  $\bar P$ described in (c), and the $(|V|+1)$-th row is the point {($\eta, \bold{x}$) = ($\frac{f_{\bold i(S )}(S )}{\alpha_{\bold i(S )}},  \sum_{i \in S } \bold {e}_i$)} given in  (a). Consider the following row operations. 
	\begin{itemize}		
	\item[Step 1:] We multiply the  $(|V|+1)$-th row by -1 to get a  row {($\frac{-f_{\bold i(S )}(S )}{\alpha_{\bold i(S )}}
		,  \sum_{i \in S } -\bold {e}_i$).  }
	\item[Step 2:] We add the new row {($\frac{-f_{\bold i(S )}(S )}{\alpha_{\bold i(S )}}
		,  \sum_{i \in S } -\bold {e}_i$)} to each of the first $| V \setminus S |$ rows. Then, we get $|V\setminus S |$ linearly independent rows, {$(
		\frac{f_{\bold i(S  \cup \{j\})}(S  \cup \{j\})}{\alpha_{\bold i(S  \cup \{j\})}}+\frac{-f_{\bold i(S )}(S )}{\alpha_{\bold i(S )}}, \bold{e}_j)$ } for all $j \in V\setminus S $.  
	\item[Step 3:] We multiply each of the $|V\setminus S |$ linearly independent rows of Step 2 by -1. We get {$(
		\frac{-f_{\bold i(S  \cup \{j\})}(S  \cup \{j\})}{\alpha_{\bold i(S  \cup \{j\})}}+\frac{f_{\bold i(S )}(S )}{\alpha_{\bold i(S )}}, -\bold{e}_j)$ } for all $j \in V\setminus S $. 
	\item[Step 4:] {Recall that from the $(| V \setminus S |+1)$-th row to the $|V|$-th row, each of the rows is represented by  {$(\frac{f_{\bold i(S  \setminus\{j\} \cup \{k_{j}\})}(S  \setminus\{j\} \cup \{k_{j}\})}{\alpha_{\bold i(S  \setminus\{j\} \cup \{k_{j}\})}}, \sum_{i \in S  \setminus\{j\}} \bold {e}_i + \bold{e}_{k_{j}})$  for a given $j \in S $. Here, for a given $j \in S $, there exists a row $ (\frac{-f_{\bold i(S  \cup \{k_j\})}(S  \cup \{k_j\})}{\alpha_{\bold i(S  \cup \{k_j\})}}+\frac{f_{\bold i(S )}(S )}{\alpha_{\bold i(S )}} , -\bold{e}_{k_j})$ from Step 3. Now for a given $j \in S $, from the $(| V \setminus S |+1)$-th row to the $V$-th row, we add  the rows $(\frac{-f_{\bold i(S )}(S )}{\alpha_{\bold i(S )}}
			,  \sum_{i \in S } -\bold {e}_i)$   and  $(\frac{-f_{\bold i(S  \cup \{k_j\})}(S  \cup \{k_j\})}{\alpha_{\bold i(S  \cup \{k_j\})}}+\frac{f_{\bold i(S )}(S )}{\alpha_{\bold i(S )}} , -\bold{e}_{k_j})$ to the row 		
		$(\frac{f_{\bold i(S  \setminus\{j\} \cup \{k_{j}\})}(S  \setminus\{j\} \cup \{k_{j}\})}{\alpha_{\bold i(S  \setminus\{j\} \cup \{k_{j}\})}}, \sum_{i \in S  \setminus\{j\}} \bold {e}_i + \bold{e}_{k_{j}})$. Then we get $|S |$ linearly independent rows, $(0,-\bold{e}_{j})$ for all $j \in S $.}}
\end{itemize}		
	Steps 1 to 4 show that the $n+1$ points described in  (a)--(c) are  affinely independent. 
\end{proof}
We provide Example \ref{ex:facet} to demonstrate Proposition \ref{prop:sub_facet}.

\begin{example}\label{ex:facet}
Suppose that we have $m = 2$ submodular functions with $\alpha_1 = \alpha_2 = 1$ and $n = 4$ elements $V= \{1,2,3,4\}$. For the case $S  = \{1,2\}$, we have two associated submodular inequalities
	\begin{align*}
	\theta_1 &\le  2 +  2x_3 + 3x_4, \text{ and}\\
	\theta_2 &\le  5 -2(1-x_1) -3(1-x_2) +x_3+4x_4,
\end{align*} 
where $f_1(S ) = 2$, $f_2(S ) = 5$,  
$\rho^{1}_1(V \setminus \{1\}) = \rho^{1}_2(V \setminus \{2\}) = 0$, $\rho^{2}_1(V \setminus \{1\}) = 2$, $\rho^{2}_2(V \setminus \{2\}) = 3$,
$\rho^{1}_3(S ) = 2$,  $\rho^{1}_4(S ) = 3$, $\rho^{2}_3(S ) = 1$, and $\rho^{2}_4(S ) = 4$. Note that the function $\bold i(S ) = \mathop{\arg\min}_{i \in [2]} \frac{f_i(S )}{\alpha_i}$ is equal to 1 since the first submodular function at $S $ attains the smallest value $f_1(S ) < f_2(S )$. Here, the submodular inequality $\eta \le  2 +  2x_3 + 3x_4$ is facet defining, because $f_1(S  \cup \{3\}) = 4$, $f_1(S  \cup \{4\}) = 5$, and $\rho^{1}_1(\{3\}) = \rho^{1}_2(\{4\}) = 0$. Condition (i) of Proposition \ref{prop:sub_facet} holds, since $\rho^{1}_1(\{3\}) = \rho^{1}_2(\{4\}) = 0$. Condition (ii) of Proposition \ref{prop:sub_facet} holds, since $f_1(S  \cup \{3\}) = f_1(S )+ \rho^{1}_3(S ) = 2+2 = 4$, $f_1(S  \cup \{4\}) = f_1(S )+ \rho^{1}_4(S ) = 2+3=5$, and $\bold i(S ) = \bold i(S  \setminus \{1\} \cup \{3\} ) = \bold i(S  \setminus \{2\} \cup \{4\} ) = 1$.

From (a)--(c) of the proof of  Proposition \ref{prop:sub_facet}, the $n+1$ affinely independent points ($\eta,x_1,x_2,x_3,x_4$) are as follows. The point (2,1,1,0,0) is based on the selection of $S $ as described in (a). From (b), there exist $|V \setminus S | = 2$ points, (4,1,1,1,0) and (5,1,1,0,1) based on the selection of $S  \cup \{3\}$ and $S  \cup \{4\}$. From (c), there exist $|S | = 2$ points (4,0,1,1,0) and (5,1,0,0,1) based on the marginal contributions $\rho^{1}_1(\{3\}) = \rho^{1}_2(\{4\}) = 0$. We demonstrate the row operation steps 1 to 4 of the proof as follows, where the final table shows that the $n+1=5$ points are affinely independent.

\(\begin{pmatrix}
	4&  1& 1& 1& 0&\\ 
	5&  1& 1& 0& 1&\\
	4&  0& 1& 1& 0&\\
	5&  1& 0& 0& 1&\\
	2&  1& 1& 0& 0&\\
\end{pmatrix} \overset{\textrm{Step 1}}{\longrightarrow} 
\begin{pmatrix}
	4&   1& 1& 1& 0&\\ 
	5&   1& 1& 0& 1&\\
	4&   0& 1& 1& 0&\\
	5&   1& 0& 0& 1&\\
	-2&  -1& -1& 0& 0&\\
\end{pmatrix} \overset{\textrm{Step 2}}{\longrightarrow} 
\begin{pmatrix}
	2&   0& 0& 1& 0&\\ 
	3&   0& 0& 0& 1&\\
	4&   0& 1& 1& 0&\\
	5&   1& 0& 0& 1&\\
	-2& -1& -1& 0& 0&\\
\end{pmatrix}\overset{\textrm{Step 3}}{\longrightarrow}
\)

\( 
\begin{pmatrix}
	-2&  0& 0& -1& 0&\\ 
	-3&  0& 0& 0& -1&\\
	4&  0& 1& 1& 0&\\
	5&   1& 0& 0& 1&\\
	-2&   -1& -1& 0& 0&\\
\end{pmatrix} \overset{\textrm{Step 4}}{\longrightarrow} 
\begin{pmatrix}
	-2& - 0& 0& -1& 0&\\ 
	-3&  - 0& 0& 0& -1&\\
	0&   -1& 0& 0& 0&\\
	0&  0& -1& 0& 0&\\
	-2&   -1& -1& 0& 0&\\
\end{pmatrix}\overset{\textrm{Final Matrix}}{\longrightarrow} 
\begin{pmatrix}
	0&  0& 0& 1& 0&\\ 
	0&  0& 0& 0& 1&\\
	0&  1& 0& 0& 0&\\
	0&  0& 1& 0& 0&\\
	1&  0& 0& 0& 0&\\
\end{pmatrix}
\)

\end{example}

We note that it may be difficult to find a submodular inequality that simultaneously meets the two conditions of Proposition \ref{prop:sub_facet}. Specifically, given a submodular inequality for $S $,  the computational effort to check whether  the two conditions hold may be close to generating all $m$ submodular inequalities corresponding to the set $S $ (not just one inequality corresponding to $\bar i$). However, we are able to derive some computational strategies based on the two conditions of Proposition \ref{prop:sub_facet}.  
 Lemma \ref{lemma:margin_equal}, Lemma \ref{prop:sub_ftn_idea2}, and Lemma \ref{prop:sub_ftn_idea2_coro} provide an important observation to this end.

{
	\begin{lemma}\label{lemma:margin_equal}
		Given an index $i \in [m]$,  $\tilde X'' \subseteq \tilde X \subseteq V$, and $\tilde S \subseteq V$, where $\tilde X$ and $\tilde S$ follow the equality $f_i(\tilde X \cup \tilde S) = f_i(\tilde S) + \sum_{j \in \tilde X} \rho^i_j(\tilde S) $, if the equality
		$\rho^i_{j}(\tilde X'' \cup \tilde S) = \rho^i_{j}(\tilde S)$ holds for $j \in \tilde X \setminus \tilde X''$, then the relation 
		\begin{equation}\label{lemma:margin_equal_2} 
			\rho^i_{j}(\tilde X'' \cup \tilde S \cup Z) = \rho^i_{j}(\tilde S \cup Z)
		\end{equation} 
		also holds for  $Z \subseteq V$.
	\end{lemma}
	\begin{proof}
		We prove the relation \eqref{lemma:margin_equal_2} by mathematical induction. Given an index $i \in [m]$ and $ \tilde S \subseteq V$, consider the base case of the induction with $\tilde X' =  \emptyset $. We have $\rho^i_{j}(\emptyset \cup \tilde S \cup Z) = \rho^i_{j}(\tilde S \cup Z)$, which trivially satisfies \eqref{lemma:margin_equal_2}, for all $j \in \tilde X \setminus \tilde X'$ and $Z \subseteq V$. Now for the  case  with  $\tilde X' = \{j_1, \dots, j_{\bar n-1}\}$ for $2\le\bar n\le |\tilde X''|$, we assume that for all $j \in \tilde X \setminus \tilde X'$ and $Z \subseteq V$, the relation 
		\begin{equation}\label{eq:inductionhyp}
			\rho^i_{j}(\tilde X' \cup \tilde S \cup Z) = \rho^i_{j}(\tilde S \cup Z)
		\end{equation} 
		holds. Now consider the  case with  $\tilde X'' = \{j_1, \dots, j_{\bar n-1},j_{\bar n}\}$. Equation \eqref{eq:inductionhyp} can be rewritten as
		\begin{equation*}
			f_i(\{j\} \cup \tilde X' \cup \tilde S \cup Z) - f_i( \tilde X' \cup \tilde S \cup Z) = f_i(\{j\} \cup \tilde S \cup Z) - f_i( \tilde S \cup Z).
		\end{equation*} 
		Note that since $j_{\bar n} \in \tilde X'' \setminus \tilde X'$, $	\rho^i_{j_{\bar n}}(\tilde X' \cup \tilde S \cup Z) = \rho^i_{j_{\bar n}}(\tilde S \cup Z)$ for $Z \subseteq V$, we can construct a new $Z' = \{j\} \cup Z \subseteq V$ and algebraically handle the above equation for all $j \in \tilde X'' \setminus \tilde X'$ as
		\begin{align*}	
			&[f_i(\{j\} \cup \tilde X' \cup \tilde S \cup Z) + \rho^i_{j_{\bar n}}(\tilde X' \cup \tilde S \cup Z')] -[f_i( \tilde X' \cup \tilde S \cup Z)+\rho^i_{j_{\bar n}}(\tilde X' \cup \tilde S \cup Z)] \\
			&= [f_i(\{j\} \cup \tilde S \cup Z)+\rho^i_{j_{\bar n}}( \tilde S \cup Z')] - [f_i( \tilde S \cup Z)+\rho^i_{j_{\bar n}}( \tilde S \cup Z)],	
		\end{align*} 
		which is equivalent to
		\begin{align*}	
			f_i(\{j\} \cup \tilde X'' \cup \tilde S \cup Z) - f_i(\tilde X'' \cup \tilde S \cup Z)
			=f_i(\{j\} \cup \tilde S \cup Z)- f_i( \tilde S \cup Z).
		\end{align*} 
		Therefore, \eqref{lemma:margin_equal_2} holds for $\tilde X''$, which completes the proof.  
	\end{proof}
}

\begin{lemma}\label{prop:sub_ftn_idea2}
	Given an index $i \in [m]$ and  $\tilde X, \tilde S \subseteq V$, if the equality
	\begin{equation}\label{prop:sub_ftn_idea2_1}
		f_i(\tilde X \cup \tilde S) = f_i(\tilde S) + \sum_{j \in \tilde X} \rho^i_j(\tilde S) 
	\end{equation} 
	holds, then the relation 
	\begin{equation}\label{prop:sub_ftn_idea2_2} 
		f_i(\tilde X \cup \tilde S \cup \{z\}) = f_i(\tilde S \cup \{z\}) + \sum_{j \in \tilde X} \rho^i_j(\tilde S \cup \{z\}) 
	\end{equation} 
	also holds for any  $z \in V$.
\end{lemma}
\begin{proof}
	We prove the relation \eqref{prop:sub_ftn_idea2_2} by mathematical induction. Given an index $i \in [m]$, $z \in V$, and $ \tilde S \subseteq V$, consider the base case of the induction with a single element $\tilde X' = \{j_1\}$, where
	\begin{equation*}
		f_i(\{j_1\} \cup \tilde S \cup \{z\})-f_i(\tilde S \cup \{z\}) = \rho^i_{j_1}(\tilde S \cup \{z\}), 
	\end{equation*}
	which follows from the definition of $\rho$. Now for the  case  with  $\tilde X'' = \{j_1, \dots, j_{\bar n-1}\}$ with $\bar n-1< |\tilde X|$ elements, we assume that under the condition 
	\begin{equation}\label{prop:case2cond}
		f_i(\tilde X'' \cup \tilde S) = f_i(\tilde S) + \sum_{j \in \tilde X''} \rho^i_j(\tilde S), 
	\end{equation}
	and the relation 
	\begin{equation}\label{prop:case2relat}
		f_i(\tilde X'' \cup \tilde S \cup \{z\}) = f_i(\tilde S \cup \{z\}) + \sum_{j \in \tilde X''} \rho^i_j(\tilde S \cup \{z\}) 
	\end{equation}  
	holds.
	
	For the  case with $\bar n$ elements $\tilde X = \{j_1, \dots, j_{\bar n-1},j_{\bar n}\}$,  we have 
	\begin{align*}	
		f_i(\tilde X'' \cup \{j_{\bar n}\} \cup \tilde S)-f_i(\tilde S ) & = f_i(\tilde X'' \cup \tilde S)+  \rho^i_{j_{\bar n}}(\tilde X'' \cup \tilde S) - f_i(\tilde S )\\
		& =  f_i(\tilde S) + \sum_{j \in \tilde X''} \rho^i_j(\tilde S) +  \rho^i_{j_{\bar n}}(\tilde X'' \cup \tilde S) - f_i(\tilde S )\\
		& = \sum_{j \in \tilde X} \rho^i_j(\tilde S),
	\end{align*} 
	where the second equality follows from  \eqref{prop:case2cond}, and the third equality  is from the condition \eqref{prop:sub_ftn_idea2_1} of the final case. Since $\sum_{j \in \tilde X''} \rho^i_j(\tilde S) +  \rho^i_{j_{\bar n}}(\tilde X'' \cup \tilde S) = \sum_{j \in \tilde X} \rho^i_j(\tilde S)$, we have $\rho^i_{j_{\bar n}}(\tilde X'' \cup \tilde S) = \rho^i_{j_{\bar n}}(\tilde S)$. {Here, the element $j_{\bar n} \in \tilde X \setminus \tilde X''$, and therefore, the relation $\rho^i_{j_{\bar n}}(\tilde S \cup \{z\})=\rho^i_{j_{\bar n}}(\tilde X'' \cup \tilde S \cup \{z\}) $ holds from Lemma \ref{lemma:margin_equal}.} We have
	\begin{align*}	
		f_i(\tilde X'' \cup \{j_{\bar n}\} \cup \tilde S \cup \{z\})-f_i(\tilde S \cup \{z\}) & = f_i(\tilde X'' \cup \tilde S \cup \{z\})+  \rho^i_{j_{\bar n}}(\tilde X'' \cup \tilde S \cup \{z\}) - f_i(\tilde S \cup \{z\})\\
		& = f_i(\tilde X'' \cup \tilde S \cup \{z\})+  \rho^i_{j_{\bar n}}(\tilde S \cup \{z\}) - f_i(\tilde S \cup \{z\}).
	\end{align*}
	
	From the above relations, since the assumption of the relation \eqref{prop:case2relat} holds, we have
	\begin{align*}	
		f_i(\tilde X \cup \tilde S \cup \{z\}) & = f_i(\tilde X'' \cup \{j_{\bar n}\} \cup \tilde S \cup \{z\}) \\
		& = f_i(\tilde X'' \cup \tilde S \cup \{z\})+  \rho^i_{j_{\bar n}}(\tilde S \cup \{z\}) \\
		& = f_i(\tilde S \cup \{z\}) + \sum_{j \in \tilde X''} \rho^i_j(\tilde S \cup \{z\}) + \rho^i_{j_{\bar n}}(\tilde S \cup \{z\}) \\
		& = f_i(\tilde S \cup \{z\}) + \sum_{j \in \tilde X} \rho^i_j(\tilde S \cup \{z\}) .
	\end{align*} 
	This completes the proof.
\end{proof}

\begin{lemma}\label{prop:sub_ftn_idea2_coro}
	Given an index $i \in [m]$ and  $\tilde X, \tilde S \subseteq V$, if the equality \eqref{prop:sub_ftn_idea2_1} of Lemma \ref{prop:sub_ftn_idea2}
	holds, then the relation 
	\begin{equation}\label{prop:sub_ftn_idea2_coro_2} 
		f_i(\tilde X \cup \tilde S \cup Z) = f_i(\tilde S \cup Z) + \sum_{j \in \tilde X} \rho^i_j(\tilde S \cup Z) 
	\end{equation} 
	also holds for $Z \subseteq V$.
\end{lemma}
\begin{proof}
	Suppose that we are given an index $i \in [m]$, $Z \subseteq V$, and $\tilde X, \tilde S$ for the equality. The following steps show that adding all elements from $Z$ to $\tilde S$ recursively does not violate the relation \eqref{prop:sub_ftn_idea2_2} of Lemma \ref{prop:sub_ftn_idea2}.
	\begin{itemize}		
		\item[Step 1:] Pick  an element $z \in Z$.
		\item[Step 2:] Since the equality \eqref{prop:sub_ftn_idea2_1} of Lemma \ref{prop:sub_ftn_idea2} holds, the relation \eqref{prop:sub_ftn_idea2_2} of Lemma \ref{prop:sub_ftn_idea2} holds.
		\item[Step 3:] Set  $\tilde S = \tilde S \cup \{z\}$ for the relation \eqref{prop:sub_ftn_idea2_2} of Lemma \ref{prop:sub_ftn_idea2}. The new $\tilde S$ satisfies the equality \eqref{prop:sub_ftn_idea2_1} of Lemma \ref{prop:sub_ftn_idea2}.
		\item[Step 3:] Let $Z = Z \setminus \{z\}$. Go to Step 1 if $Z \neq \emptyset$; otherwise, stop.
	\end{itemize}	
	This completes the proof.
\end{proof}

Using this lemma, we provide a proposition that informs a useful computational strategy to select a more compact set of sufficient submodular inequalities.  
We separate the set $S $ into two disjoint subsets. From the first subset, we derive a new subset  of elements, which is based on  condition (i) of Proposition \ref{prop:sub_facet}. Then, we consider a  union of the second subset with the new subset and make sure that the submodular inequality \eqref{eq:submod-eta} associated with this particular union of subsets does not violate Proposition \ref{prop:MIP_valid}.

\begin{proposition}\label{prop:inter_cut}
	Given a set $\bar X \subseteq V$ and an index $i \in [m]$, we define two associated subsets $\tilde X_i \subseteq \bar X$ and
\begin{equation}\label{prop:inter_cut1}
	 \mathcal{S}(i,\tilde X_i) = \{ j \in V\setminus \tilde X_i:  \exists k \in \tilde X_i \mbox{ with } \rho^{i}_j(\{k\}) = 0 \}.
\end{equation}
 If the condition 
\begin{equation}\label{prop:inter_cut2}
	f_{\bold i(\bar X)}(\tilde X_{\bold i(\bar X)})=f_{\bold i(\bar X)}(\mathcal{S}(\bold i(\bar X),\tilde X_{\bold i(\bar X)}))+ \sum_{j \in \tilde X_{\bold i(\bar X)}} \rho^{\bold i(\bar X)}_j(\mathcal{S}(\bold i(\bar X),\tilde X_{\bold i(\bar X)}))
\end{equation}
	holds for all $\bar X\subseteq V$, then  using  the mixed-integer set given   by 
	\begin{align*}
		\mathcal F'= 	&\{(\eta,\bold{x})\in \R \times \mathbb B^{n}:  \\ 
		& \eta \le \frac{1}{\alpha_{\bold i(\bar X)}} ( f_{\bold i(\bar X)}( \mathcal{S}(\bold i(\bar X),\tilde X_{\bold i(\bar X)}) \cup \bar X \setminus \tilde X_{\bold i(\bar X)}) - 
		\sum_{j\in \mathcal{S}(\bold i(\bar X),\tilde X_{\bold i(\bar X)})\cup \bar X \setminus \tilde X_{\bold i(\bar X)}}  \rho^{\bold i(\bar X)}_j(V\setminus\{j\})(1-x_j)+\\
		&\sum_{j\in V\setminus \{\mathcal{S}(\bold i(\bar X),\tilde X_{\bold i(\bar X)})\cup \bar X \setminus \tilde X_{\bold i(\bar X)} \}} \rho^{\bold i(\bar X)}_j(\mathcal{S}(\bold i(\bar X),\tilde X_{\bold i(\bar X)})\cup \bar X \setminus \tilde X_{\bold i(\bar X)}) x_j), \forall \bar  X\subseteq V\},
	\end{align*}
  to define  the set $\mathcal{C}$ in Formulation \eqref{RMP}  provides an optimal solution to Problem \eqref{RSM}.
\end{proposition} 	
\begin{proof}
Given $\bar X \subseteq V$ and the associated $\tilde X_{\bold i(\bar X)} \subseteq \bar X$, we have 
\begin{subequations}\label{prop:inter_cutproof}
\begin{align} 	
	&f_{\bold i(\bar X)}(\bar X)-f_{\bold i(\bar X)}(\mathcal{S}(\bold i(\bar X),\tilde X_{\bold i(\bar X)}) \cup \bar  X \setminus \tilde X_{\bold i(\bar X)} ) \notag\\ & = f_{\bold i(\bar X)}(\tilde X_{\bold i(\bar X)} \cup \mathcal{S}(\bold i(\bar X),\tilde X_{\bold i(\bar X)}) \cup \bar X \setminus \tilde X_{\bold i(\bar X)})-f_{\bold i(\bar X)}(\mathcal{S}(\bold i(\bar X),\tilde X_{\bold i(\bar X)}) \cup \bar  X \setminus \tilde X_{\bold i(\bar X)} )\label{prop:inter_cutproof1}\\
	&= \sum_{j \in \tilde X_{\bold i(\bar X)}} \rho^{\bold i(\bar X)}_j(\mathcal{S}(\bold i(\bar X),\tilde X_{\bold i(\bar X)})\cup \bar  X \setminus \tilde X_{\bold i(\bar X)}).\label{prop:inter_cutproof2}
\end{align} 
\end{subequations}	Equality \eqref{prop:inter_cutproof1} follows from $\tilde X_{\bold i(\bar X)} \subseteq \bar X$ and $f_{\bold i(\bar X)}(\tilde X_{\bold i(\bar X)}) = f_{\bold i(\bar X)}(\tilde X_{\bold i(\bar X)} \cup \mathcal{S}(\bold i(\bar X),\tilde X_{\bold i(\bar X)}))$ since for all $j \in \mathcal{S}(\bold i(\bar X),\tilde X_{\bold i(\bar X)})$, there exists $k \in \tilde X_{\bold i(\bar X)}$ such that $\rho^{\bold i(\bar X)}_j(\{k\}) = 0$ shown in the definition of \eqref{prop:inter_cut1}. Equality \eqref{prop:inter_cutproof2} follows from {condition \eqref{prop:inter_cut2} and Lemma \ref{prop:sub_ftn_idea2_coro} as follows. Suppose that $X' = \bar X \setminus \tilde X_{\bold i(\bar X)}$. Equality \eqref{prop:inter_cutproof2} provides
\begin{align*} 	
	f_{\bold i(\bar X)}(\mathcal{S}(\bold i(\bar X),\tilde X_{\bold i(\bar X)})\cup \tilde X_{\bold i(\bar X)} \cup X' )& =f_{\bold i(\bar X)}(\bar X) \\ 
	&= \sum_{j \in \tilde X_{\bold i(\bar X)}} \rho^{\bold i(\bar X)}_j(\mathcal{S}(\bold i(\bar X),\tilde X_{\bold i(\bar X)})\cup X' )+f_{\bold i(\bar X)}(\mathcal{S}(\bold i(\bar X),\tilde X_{\bold i(\bar X)}) \cup \bar  X \setminus \tilde X_{\bold i(\bar X)} )\\
	&= \sum_{j \in \tilde X_{\bold i(\bar X)}} \rho^{\bold i(\bar X)}_j(\mathcal{S}(\bold i(\bar X),\tilde X_{\bold i(\bar X)})\cup  X')+f_{\bold i(\bar X)}(\mathcal{S}(\bold i(\bar X),\tilde X_{\bold i(\bar X)}) \cup  X' ),
\end{align*} 
where we note that $\mathcal{S}(\bold i(\bar X),\tilde X_{\bold i(\bar X)})$ is $\tilde S$ of Lemma \ref{prop:sub_ftn_idea2_coro}, $\tilde X_{\bold i(\bar X)}$ is $\tilde X$ of Lemma \ref{prop:sub_ftn_idea2_coro}, and $X'$ is $Z$ of Lemma \ref{prop:sub_ftn_idea2_coro}.
}
Consider the given $\bar X \subseteq V$ and $\tilde X_{\bold i(\bar X)}$ with the subset $\mathcal{S}(\bold i(\bar X),\tilde X_{\bold i(\bar X)})$ that follows  \eqref{prop:inter_cut1} and \eqref{prop:inter_cut2}. Then
\begin{equation*}
\eta \le \frac{1}{\alpha_{\bold i(\bar X)}} ( f_{\bold i(\bar X)}( \mathcal{S}(\bold i(\bar X),\tilde X_{\bold i(\bar X)}) \cup \bar X \setminus \tilde X_{\bold i(\bar X)})+\sum_{j\in \tilde X_{\bold i(\bar X)}} \rho^{\bold i(\bar X)}_j(\mathcal{S}(\bold i(\bar X),\tilde X_{\bold i(\bar X)})\cup \bar X \setminus \tilde X_{\bold i(\bar X)}) ) = \frac{f_{\bold i(\bar X)}(\bar X) }{\alpha_{\bold i(\bar X)}},
\end{equation*} where the equality follows from \eqref{prop:inter_cutproof} and 
\eqref{prop:inter_cut1} with $\rho^{\bold i(\bar X)}_j(V\setminus\{j\}) = 0$ for all $j \in \mathcal{S}(\bold i(\bar X),\tilde X_{\bold i(\bar X)})$. Finally, {Formulation \eqref{RMP} with the mixed-integer set $\mathcal F'$ derived from Formulation \eqref{RSM_alt} provides  
\begin{align*}
	\eta & \le \frac{\theta_{\bold i(\bar X)}}{\alpha_{\bold i(\bar X)}}  \le   \frac{f_{\bold i(\bar X)}(\bar X)}{\alpha_{\bold i(\bar X)}}\\
	&\le \frac{f_{\bold i(\bar X)}( \mathcal{S}(\bold i(\bar X), \tilde X_{\bold i(\bar X)}) \cup \bar  X \setminus \tilde X_{\bold i(\bar X)}) - \sum_{j\in \mathcal{S}(\bold i(\bar X),\tilde X_{\bold i(\bar X)})\cup \bar X \setminus \tilde X_{\bold i(\bar X)}}  \rho^{\bold i(\bar X)}_j(V\setminus\{j\})(1-x_j) }{\alpha_{\bold i(\bar X)}}\\ 
	&+ \frac{\sum_{j\in V\setminus \{\mathcal{S}(\bold i(\bar X),\tilde X_{\bold i(\bar X)})\cup \bar  X \setminus \tilde X_{\bold i(\bar X)} \}} \rho^{\bold i(\bar X)}_j(\mathcal{S}(\bold i(\bar X),\tilde X_{\bold i(\bar X)})\cup \bar X \setminus \tilde X_{\bold i(\bar X)}) x_j }{\alpha_{\bold i(\bar X)}}. 
\end{align*} }
Following the end of Proposition \ref{prop:MIP_valid}, this completes the proof.
\end{proof} 

In Proposition \ref{prop:inter_cut}, from  condition (i) of Proposition \ref{prop:sub_facet}, we define a set  $\mathcal{S}(i,\tilde X_i) = \{ j \in V\setminus \tilde X_i: \exists k \in \tilde X_i  \mbox{ with }\rho^{i}_j(\{k\}) = 0\}$ based on an index $i \in [m]$ and a subset $\tilde X_i \subseteq \bar X \subseteq V$, where an index $j \in \mathcal{S}(i,\tilde X_i)$ has at least one associated index $k_j \in \tilde X_i$  such that $\rho^{\bar i}_j(\{k_j\}) = 0$. Then, if  condition \eqref{prop:inter_cut2} is satisfied, we show that with $i = \bold i(\bar X)$ and $S = \mathcal{S}(\bold i(\bar X),\tilde X_{\bold i(\bar X)}) \cup \bar X \setminus \tilde X_{\bold i(\bar X)}$, the associated submodular inequality \eqref{eq:submod-eta} provides an upper bound $   \frac{f_{\bold i(\bar X)}(\bar X)}{\alpha_{\bold i(\bar X)}}$ of the RSM \eqref{RSM_alt} for a solution $\bold{\bar x} \in \mathcal{X}$. The verification of the upper bound for a solution is necessary to establish that it suffices to consider the set $\mathcal F'$ in defining set $\mathcal C$. This further enhances the computational efficiency, as we will show in our computational study.

We now consider condition (ii) of Proposition \ref{prop:sub_facet}. Although finding a facet-defining submodular inequality is challenging, we give  a sequence of two propositions showing when, under certain conditions, a submodular inequality \eqref{eq:submod-eta} is redundant and when it is a facet of  conv$(\mathcal F)$. We first show that given $\bar X \subseteq V$, some submodular inequalities based on $\bar X$ may be redundant (i.e., dominated) in  RMP \eqref{RMP}.

\begin{proposition}\label{coro:sub_redundant}
	Given $\bar X \subseteq V$ and $\bar i,\bar i' \in [m]$ and $\bar i \neq \bar i'$,  if  $\frac{f_{\bar i}(\bar X) }{\alpha_{\bar i}} \le \frac{f_{\bar i'}(\bar X)}{\alpha_{\bar i'}}$ , $\frac{-\rho^{\bar i}_j(V\setminus\{j\})}{\alpha_{\bar i}} \le \frac{-\rho^{\bar i'}_j(V\setminus\{j\})}{\alpha_{\bar i'}}$  for all $j\in \bar X$, and
		$\frac{\rho^{\bar i}_j(\bar X) }{\alpha_{\bar i}} \le \frac{\rho^{\bar i'}_j(\bar X)}{\alpha_{\bar i'}}$ for all $j\in V\setminus \bar X$,
	then  inequality \eqref{eq:submod-eta}  with $\bar X \subseteq V$ and  $i = \bar i'$,
	\begin{equation*}
		\eta \le \frac{1}{\alpha_{\bar i'}}
		(f_{\bar i'}(\bar X) - \sum_{j\in \bar X}  \rho^{\bar i'}_j(V\setminus\{j\})(1-x_j)+\sum_{j\in V\setminus \bar X} \rho^{\bar i'}_j(\bar X) x_j),
	\end{equation*} 
	is redundant in RMP \eqref{RMP}.
\end{proposition} 
\begin{proof} 
	We follow the proof of Proposition \ref{prop:MIP_valid} and the relations of Proposition \ref{coro:sub_redundant}, and obtain
	\begin{align*}
		\eta & \le \frac{\theta_{\bold i(\bar X)}}{\alpha_{\bold i(\bar X)}}  \le   \frac{f_{\bold i(\bar X)}(\bar X)}{\alpha_{\bold i(\bar X)}}\\
		& \le \frac{ f_{\bar i}(\bar X) - \sum_{j\in \bar X}  \rho^{\bar i}_j(V\setminus\{j\})(1-x_j)+\sum_{j\in V\setminus \bar X} \rho^{\bar i}_j(\bar X) x_j}{\alpha_{\bar i}}\\
		& \le \frac{ f_{\bar i'}(\bar X) - \sum_{j\in \bar X}  \rho^{\bar i'}_j(V\setminus\{j\})(1-x_j)+\sum_{j\in V\setminus \bar X} \rho^{\bar i'}_j(\bar X) x_j}{\alpha_{\bar i'}},	
	\end{align*}  
	where inequality \eqref{eq:submod-eta}  with $\bar X \subseteq V$ and  $i = \bar i$ provides a better upper bound compared to the submodular inequality \eqref{eq:submod-eta}  with $\bar X \subseteq V$ and  $i = \bar i'$. This completes the proof.
\end{proof}

\begin{example}\label{ex:redundant}
	Suppose that we have $m = 3$ submodular functions with $\alpha_1 = \alpha_2 = \alpha_3= 1$ and $n = 3$ elements $V= \{1,2,3,4\}$. For  $\bar X = \{1,2\}$, we have three associated submodular inequlities 
	\begin{align*}
		\eta &\le 3 + 2x_3 + 3x_4, \text{ and}\\
		\eta &\le 2 + 3x_3 + 4x_4, \text{ and}\\
		\eta &\le 5 + 3x_3 + 5x_4 \}.
	\end{align*} 
	The third inequality  is redundant for  RMP \eqref{RMP} since $3 + 2x_3 + 3x_4 \le 5 + 3x_3 + 5x_4$ and $2+ 3x_3 + 4x_4 \le 5 + 3x_3 + 5x_4$  from Proposition \ref{coro:sub_redundant}.
\end{example}

Proposition \ref{coro:sub_redundant} shows that if a submodular inequality's right-hand side (RHS) and coefficients are all greater than those of  another submodular inequality, the former is redundant for  RMP \eqref{RMP}. Based on Propositions  \ref{prop:sub_facet} and \ref{coro:sub_redundant}, we also give a corollary that, under certain conditions,  a given set of submodular inequalities is  a {facet} for  RMP \eqref{RMP}. Given a subset $S \subseteq V$ and $I \subseteq [m]$, we define a {mixed-integer set of the set of submodular inequalities} as $C (S, I) = \{
(\eta,\bold{x})\in \R \times \mathbb B^{n}:\eta \le \frac{1}{\alpha_{i}}( f_{i}(S) - \sum_{j\in S}  \rho^{i}_j(V\setminus\{j\})(1-x_j)+\sum_{j\in V\setminus S} \rho^{i}_j(S) x_j), \forall i \in I \}$.

\begin{corollary}\label{coro:sub_cuts_facet}
	Given  $\bar X \subseteq V$ and $I \subseteq [m]$, each submodular inequality defining the set $C (\bar X, I) = \{
	(\eta,\bold{x})\in \R \times \mathbb B^{n}:\eta \le \frac{1}{\alpha_{i}} (f_{i}(\bar X) - \sum_{j\in \bar X}  \rho^{i}_j(V\setminus\{j\})(1-x_j)+\sum_{j\in V\setminus \bar X} \rho^{i}_j(\bar X) x_j), \forall i \in I \}$ is a facet of  conv($\mathcal F$) if the following conditions hold 
	\begin{itemize}		
	\item[(i)] {for all $j \in \bar X$ and $i \in I$, there  exists at least an element $k_j \in V \setminus \bar X $ such that $\rho^{i}_j(\{k_j\}) = 0$  and  $\frac{f_{i}(\bar X)}{\alpha_{i}} = \frac{f_{\bold i(\bar X)}(\bar X)}{\alpha_{\bold i(\bar X)}}
			= \frac{f_{\bold i(\bar X \setminus \{j\} \cup \{k_j\})}(\bar X \setminus \{j\} \cup \{k_j\})}{\alpha_{\bold i(\bar X \setminus \{j\} \cup \{k_j\})}}
			= \frac{f_{\bold i(\bar X \cup \{k_j\})}(\bar X \cup \{k_j\})}{\alpha_{\bold i(\bar X \cup \{k_j\})}}$,}	
	\item[(ii)] {for all $i \in  I$, we have  $
	\frac{f_{\bold i(\bar X)}(\bar X)}{\alpha_{\bold i(\bar X)}}	 = \frac{f_{i}(\bar X)}{\alpha_i} $, } 
	\item[(iii)] {for any $\bar i' \in [m]\setminus  I$, given an inequality $(\eta, \bold x) \in C (\bar X, [m]\setminus  I)$, there must exist an index $\bar i \in  I$ for another inequality $(\eta, \bold x) \in C (\bar X,  I)$ such that  the relations of Proposition \ref{coro:sub_redundant} hold.}
\end{itemize}
\end{corollary}
\begin{proof}
	 Condition (i) delineates that all submodular inequalities in $C (\bar X, I)$ satisfy  condition (i) of Proposition \ref{prop:sub_facet} for all $\bar i \in I$. Conditions (ii) and (iii) imply that given  $\bar X \subseteq V$ and $j \in V \setminus \bar X$, we have $\min_{i \in I} \{\frac{f_{i}(\bar X)+\rho^{i}_j(\bar X)}{\alpha_i}\}= \frac{f_{\bold i(\bar X \cup \{j\})}(\bar X \cup \{j\})}{\alpha_{\bold i(\bar X \cup \{j\}})}$ since $C (\bar X, [m]\setminus  I)$ {includes redundant inequalities} (from Proposition
	\ref{coro:sub_redundant}). Then $n+1$ affinely independent points defined in (a)-(c) in the proof of Proposition \ref{prop:sub_facet} satisfy conditions (i)--(iii) of this corollary.
\end{proof}

Corollary \ref{coro:sub_cuts_facet} shows that given  $\bar X \subseteq V$, if the RHS of $\bar X$ satisfies  condition (i) of Proposition \ref{prop:sub_facet}, it may not be necessary to include all submodular inequalities for all $i \in [m]$. Next, we derive the following corollary directly from Corollary \ref{coro:sub_cuts_facet}.
\begin{corollary}\label{coro:empty_set_cut}
	The  inequalities defined by $C (\emptyset, [m])$ are facets of conv($\mathcal F$).
\end{corollary}
\begin{proof}	
Since $\bar X= \emptyset$, for any $k \in V$, it follows from  condition (i) of Corollary \ref{coro:sub_cuts_facet} that $\rho^{\bar i}_j(\{k\}) = 0$ for all $j \in \bar X$ and $i \in [m]$. Furthermore,  conditions (ii) and (iii) hold, because $I = [m]$ and $f_i(\emptyset) = 0$ for all $i \in I$.

\end{proof}

\begin{example}\label{ex:facet_set}
	Suppose that we have $m = 3$ submodular functions with $\alpha_1 = \alpha_2 = \alpha_3= 1$ and $n = 3$ elements $V= \{1,2,3\}$. For the case $\bar X = \emptyset$, we have three associated submodular inequalities $C (\emptyset, [m]) = \{$
	\begin{align*}
		\eta &\le 0+  2x_1 + 2x_2 + 3x_3, \text{ and}\\
		\eta &\le 0+ x_1 + 3x_2 + 4x_3, \text{ and}\\
		\eta &\le 0+ 3x_1 + 3x_2 + x_3 \}.
	\end{align*} 
	For the point ($x_1,x_2,x_3$) = (1,0,0), the second inequality provides an upper bound equal to 1 for the variables $\eta$ and $\theta_2$. 	For the point ($x_1,x_2,x_3$) = (0,1,0), the first inequality provides an upper bound equal to 2 for  $\eta$. For the point ($x_1,x_2,x_3$) = (0,0,1), the third inequality provides an upper bound equal to 1 for $\eta$. The $n+1$ affinely independent points ($\eta,x_1,x_2,x_3$) are (1,1,0,0), (2,0,1,0), (1,0,0,1), and (0,0,0,0).
\end{example}

\subsection{An Analysis of a Special Case of  RSM}\label{sec:AnaRSM3}

At the end of Section \ref{problem}, we highlighted the difficulty of solving  Problem \eqref{RSM3}. That is, to get $m$ values $f_i(\bold{x}_i^*)$ for all $i \in [m]$, we have to solve $m$ $\mathcal{NP}-$hard problems \eqref{SubMax}. {Let $F_i$ be a mixed-integer set defined by the set of submodular inequalities for each $i \in [m]$, i.e., $F_i=\{(\theta_i,\bold{x})\in \R \times \mathbb B^{n}:\theta_i \le f_i(S) - \sum_{j\in S}  \rho^i_j(V\setminus\{j\})(1-x_j)+\sum_{j\in V\setminus S} \rho^i_j(S) x_j , \forall S\subseteq V\}$. Recall that  $\bold {x_i^*}$
is the optimal solution to the $i$-th submodular maximization problem 
\eqref{SubMax} and $f_i(\bold{x}_i^*) = \max\{\theta_i: (\theta_i, \bold{x}) \in F_i, \bold{x} \in \mathcal{X}\cap \mathbb B^{n}, \theta_i \in \mathbb R\}$ for all $i \in [m]$.}   Let $LB$ and $UB$ be  lower and upper bounds of the optimal value of Problem \eqref{RSM3}, respectively, i.e., $LB \le \max_{\bold{x} \in \mathcal{X}\cap \mathbb B^{n}}  \min_{i \in [m]} \frac{f_i(\bold{x})}{f_i(\bold{x}_i^*)} \le UB$. It may appear that, without solving the $m$ problems, we cannot solve Problem \eqref{RSM3} or even find an optimality gap $\frac{UB-LB}{UB}$. We show how we can overcome this difficulty based on the following proposition.

{
\begin{proposition}\label{Bound_RSM2}
	Let $lb_i$ and  $ub_i$ be  lower and upper bounds of the optimal value of the $i$-th Problem \eqref{SubMax} for all $i \in [m]$, i.e., $lb_i \le f_i (\bold {x_i^*}) \le ub_i$. Let $\bar \eta_{relax} = \max\{\eta: \eta \le \frac{\theta_i}{lb_i},\  i \in [m],\ (\theta_i, \bold{x}) \in  F_i, {i \in [m]}, \ \bold{x} \in \mathcal{X}, \ \eta \in \mathbb R, \ \boldsymbol{\theta} \in \mathbb R^{m}\}$ be the  objective value of the relaxation of RSM \eqref{RSM}. For a given  $\bold{\bar x} \in \mathcal{X}\cap \mathbb B^{n}$, we have 
\begin{equation}\label{eq:Bound_RSM2}
	 	\min_{i \in [m]} \frac{f_i(\bold{\bar x})}{ub_i} \le \max_{\bold{x} \in \mathcal{X}\cap \mathbb B^{n}}  \min_{i \in [m]} \frac{f_i(\bold{x})}{f_i(\bold{x}_i^*)} 
	 	\le \bar \eta_{relax}.
\end{equation} 	
\end{proposition}
}
\begin{proof}
We start by showing that $\max_{\bold{x} \in \mathcal{X}\cap \mathbb B^{n}}  \min_{i \in [m]} \frac{f_i(\bold{x})}{f_i(\bold{x}_i^*)}
	\le \bar \eta_{relax}$. From Formulation \eqref{RSM_alt}, we have $\max \{ \eta : \eta \le \frac{f_i(\bold{x})}{f_i(\bold{x}_i^*)} \ \forall i \in [m], \bold{x} \in \mathcal{X}\cap \mathbb B^{n}, \eta \in \mathbb R \} = \max_{\bold{x} \in \mathcal{X}\cap \mathbb B^{n}}  \min_{i \in [m]} \frac{f_i(\bold{x})}{f_i(\bold{x}_i^*)}$. 
Since  $lb_i \le f_i (\bold {x_i^*})$ for all $i \in [m]$, the constraint $\eta \le \frac{f_i(\bold{x})}{lb_i}$ is a relaxation of  $\eta \le \frac{f_i(\bold{x})}{f_i(\bold{x}_i^*)}$ for Formulation \eqref{RSM_alt}. We have the following inequality
\begin{equation*}
\max \{ \eta : \eta \le \frac{f_i(\bold{x})}{f_i(\bold{x}_i^*)},   i \in [m], \bold{x} \in \mathcal{X}\cap \mathbb B^{n}, \eta \in \mathbb R \} 
\le \max \{ \eta : \eta \le \frac{f_i(\bold{x})}{lb_i}, i \in [m], \bold{x} \in \mathcal{X}\cap \mathbb B^{n}, \eta \in \mathbb R \}.
\end{equation*}	
Furthermore, the objective value $\bar \eta_{relax}$ is obtained from the relaxation of $\bold{x} \in \mathcal{X}\cap \mathbb B^{n}$ as $\bold{x} \in \mathcal{X}$. Therefore, we conclude that 
	\begin{align*}
	\max_{\bold{x} \in \mathcal{X}\cap \mathbb B^{n}}  \min_{i \in [m]} \frac{f_i(\bold{x})}{f_i(\bold{x}_i^*)}  
	& = \max \{ \eta : \eta \le \frac{f_i(\bold{x})}{f_i(\bold{x}_i^*)}, i \in [m], \bold{x} \in \mathcal{X}\cap \mathbb B^{n}, \eta \in \mathbb R \} \\
	& \le \max \{ \eta : \eta \le \frac{f_i(\bold{x})}{lb_i}, i \in [m], \bold{x} \in \mathcal{X}, \eta \in \mathbb R \} \\
	& \le \max\{\eta: \eta \le \frac{\theta_i}{lb_i}, i \in [m],\ (\theta_i, \bold{x}) \in  F_i, {i \in [m]}, \ \bold{x} \in \mathcal{X}, \ \eta \in \mathbb R, \ \boldsymbol{\theta} \in \mathbb R^{m}\} \\
	& = \bar \eta_{relax}.
\end{align*} 

Next, we show the second part of inequality \eqref{eq:Bound_RSM2}, $	 	\min_{i \in [m]} \frac{f_i(\bold{\bar x})}{ub_i} \le \max_{\bold{x} \in \mathcal{X}\cap \mathbb B^{n}}  \min_{i \in [m]} \frac{f_i(\bold{x})}{f_i(\bold{x}_i^*)} $. Since  $f_i (\bold {x_i^*}) \le ub_i$ for all $i \in [m]$, we have $\eta \le \frac{f_i(\bold{x})}{ub_i} \le \frac{f_i(\bold{x})}{f_i(\bold{x}_i^*)}$.  Thus, 
\begin{equation*}
	\max \{ \eta : \eta \le \frac{f_i(\bold{x})}{ub_i},  i \in [m], \bold{x} \in \mathcal{X}\cap \mathbb B^{n}, \eta \in \mathbb R \}  \le  \max \{ \eta : \eta \le \frac{f_i(\bold{x})}{f_i(\bold{x}_i^*)},  i \in [m], \bold{x} \in \mathcal{X}\cap \mathbb B^{n}, \eta \in \mathbb R \}.
\end{equation*}	
In addition, the solution $\bold{\bar x} \in \mathcal{X}\cap \mathbb B^{n}$ satisfies 
\begin{equation*}
	  \max \{ \eta : \eta \le \frac{f_i(\bold{\bar x})}{ub_i}, i \in [m], \eta \in \mathbb R \} \le \max \{ \eta : \eta \le \frac{f_i(\bold{x})}{ub_i}, i \in [m], \bold{x} \in \mathcal{X}\cap \mathbb B^{n}, \eta \in \mathbb R \}.
\end{equation*}		 
From the  above relations, we conclude 
	\begin{align*}
	\min_{i \in [m]} \frac{f_i(\bold{\bar x})}{ub_i} & = 
	\max \{ \eta : \eta \le \frac{f_i(\bold{\bar x})}{ub_i}, i \in [m], \eta \in \mathbb R \} \\	
	& \le \max \{ \eta : \eta \le \frac{f_i(\bold{x})}{f_i(\bold{x}_i^*)}, i \in [m], \bold{x} \in \mathcal{X}\cap \mathbb B^{n}, \eta \in \mathbb R \} \\
	& = \max_{\bold{x} \in \mathcal{X}\cap \mathbb B^{n}}  \min_{i \in [m]} \frac{f_i(\bold{x})}{f_i(\bold{x}_i^*)} .
\end{align*} 
This completes the proof.
\end{proof}

Here, we also make an observation that we can solve Problem \eqref{RSM3} without exactly solving $m$ submodular maximization problems, by instead solving  Problem \eqref{RSM} with a particular choice of $\boldsymbol{\alpha}$, such that $lb_i \le \alpha_i \le ub_i$ for all $i \in [m]$, under certain conditions.

\begin{proposition}\label{prop:no_needs_all_optimal}
	Let $\bold{\bar x'}$ be an optimal solution of $\max_{\bold{x} \in \mathcal{X}\cap \mathbb B^{n}}  \min_{i \in [m]} \frac{f_i(\bold{x})}{lb_i}$ and $\bold i(\bar X') = \mathop{\arg\min}_{i \in [m]} \frac{f_i(\bar X')}{lb_i}$. If $lb_{\bold i(\bar X')} \ge ub_i$ for all $i \in [m] \setminus \{\bold i(\bar X')\}$, then $\bold{\bar x'}$ is an optimal solution of Problem \eqref{RSM3}.
\end{proposition}
\begin{proof}		
	Since $f_{\bold i(\bar X')}(\bold{x}_{\bold i(\bar X')}^*) \ge lb_{\bold i(\bar X')} \ge ub_i$ for all $i \in [m] \setminus \{\bold i(\bar X')\}$ and $\max_{\bold{x} \in \mathcal{X}\cap \mathbb B^{n}}  \min_{i \in [m]} \frac{f_i(\bold{x})}{lb_i} = \min_{i \in [m]} \frac{f_i(\bold{\bar x'})}{lb_i}$, we have the following relation
\begin{equation*}
\frac{f_{\bold i(\bar X')}(\bold{\bar x'})}{f_{\bold i(\bar X')}(\bold{x}_{\bold i(\bar X')}^*)} \le \frac{ f_{\bold i(\bar X') }(\bold{\bar x'})}{lb_{\bold i(\bar X')}} \le  \frac{f_i(\bold{\bar x'})}{ub_i} \le \frac{f_i(\bold{\bar x'})}{f_i(\bold{x}_i^*)} \le \frac{f_i(\bold{\bar x'})}{lb_i},  i \in [m] \setminus \{\bold i(\bar X')\}.
\end{equation*} 
Consequently, the formulation $\max_{\bold{x} \in \mathcal{X}\cap \mathbb B^{n}}  \min_{i \in [m]} \frac{f_i(\bold{x})}{lb_i}$ has the following relation with several optimization problems 
	\begin{align*}	
	\max_{\bold{x} \in \mathcal{X}\cap \mathbb B^{n}}  \min_{i \in [m]} \frac{f_i(\bold{x})}{lb_i} & = \min_{i \in [m]} \frac{f_i(\bold{\bar x'})}{lb_i} \\
	& = \max \{ \eta : \eta \le \frac{f_i(\bold{\bar x'})}{lb_i}, i \in [m], \eta \in \mathbb R \} \\	
	& =	\max \{ \eta : \eta \le \frac{f_{\bold i(\bar X')}(\bold{\bar x'})}{lb_{\bold i(\bar X')}},
	\eta \le \frac{f_i(\bold{\bar x'})}{lb_i}, i \in [m] \setminus \{\bold i(\bar X')\}, \eta \in \mathbb R \}\\
	& =	\max \{ \eta : \eta \le \frac{f_{\bold i(\bar X')}(\bold{\bar x'})}{lb_{\bold i(\bar X')}},
	\eta \le \frac{f_i(\bold{\bar x'})}{ub_i}, i \in [m] \setminus \{\bold i(\bar X')\}, \eta \in \mathbb R\}\\
	& =	\max \{ \eta : \eta \le \frac{f_{\bold i(\bar X')}(\bold{\bar x'})}{lb_{\bold i(\bar X')}},
\eta \le \frac{f_i(\bold{\bar x'})}{f_i(\bold{x}_i^*)}, i \in [m] \setminus \{\bold i(\bar X')\}, \eta \in \mathbb R\}\\	
	&=  \max \{ \eta : \eta \le \frac{f_{\bold i(\bar X')}(\bold{ x})}{lb_{\bold i(\bar X')}}, \eta \le \frac{f_i(\bold{x})}{f_i(\bold{x}_i^*)}, i \in [m]\setminus \{\bold i(\bar X')\}, \bold{x} \in \mathcal{X}\cap \mathbb B^{n}, \eta \in \mathbb R \}. 
\end{align*} 
{From the above relations, since $ \frac{f_{\bold i(\bar X')}(\bold{ x})}{lb_{\bold i(\bar X')}} \le \frac{f_i(\bold{x})}{f_i(\bold{x}_i^*)}$ for all $i \in [m]$, the inequalities $\{ \eta \le \frac{f_i(\bold{x})}{f_i(\bold{x}_i^*)} , \ \forall i \in [m]\setminus \{\bold i(\bar X')\} \}$ are redundant while solving $\max \{ \eta : \eta \le \frac{f_i(\bold{\bar x'})}{lb_i} \ \forall i \in [m], \eta \in \mathbb R \}$. Therefore,}
\begin{align*}	
	\max_{\bold{x} \in \mathcal{X}\cap \mathbb B^{n}}  \min_{i \in [m]} \frac{f_i(\bold{x})}{lb_i} & = \min_{i \in [m]} \frac{f_i(\bold{\bar x'})}{lb_i} \\
	&=  \max \{ \eta : \eta \le \frac{f_{\bold i(\bar X')}(\bold{ x})}{lb_{\bold i(\bar X')}}, \bold{x} \in \mathcal{X}\cap \mathbb B^{n}, \eta \in \mathbb R \}. \\
\end{align*} 
{Since  $f_{\bold i(\bar X')}(\bold{x}_{\bold i(\bar X')}^*) \ge lb_{\bold i(\bar X')} \ge ub_i \ge f_i(\bold{x}_i^*)$ for all $i \in [m] \setminus \{\bold i(\bar X')\}$, the inequalities $\{\eta \le \frac{f_i(\bold{x})}{f_i(\bold{x}_i^*)} \ \forall i \in [m]\setminus \{\bold i(\bar X')\} \}$ are  redundant	for $\max \{ \eta :  \eta \le \frac{f_i(\bold{x})}{f_i(\bold{x}_i^*)} \ \forall i \in [m], \bold{x} \in \mathcal{X}\cap \mathbb B^{n}, \eta \in \mathbb R \}$. Therefore,}
\begin{align*}	
	\max_{\bold{x} \in \mathcal{X}\cap \mathbb B^{n}}  \min_{i \in [m]} \frac{f_i(\bold{x})}{f_i(\bold{x}_i^*)} 
&=  \max \{ \eta : \eta \le \frac{f_{\bold i(\bar X')}(\bold{ x})}{f_{\bold i(\bar X')}(\bold{x}_{\bold i(\bar X')}^*)}, \eta \le \frac{f_i(\bold{x})}{f_i(\bold{x}_i^*)}, i \in [m]\setminus \{\bold i(\bar X')\}, \bold{x} \in \mathcal{X}\cap \mathbb B^{n}, \eta \in \mathbb R \}\\ 	
	&=  \max \{ \eta : \eta \le \frac{f_{\bold i(\bar X')}(\bold{ x})}{f_{\bold i(\bar X')}(\bold{x}_{\bold i(\bar X')}^*)}, \bold{x} \in \mathcal{X}\cap \mathbb B^{n}, \eta \in \mathbb R \}.
\end{align*} 
From the above relations, we observe that  $lb_{\bold i(\bar X')}$ and $f_{\bold i(\bar X')}(\bold{x}_{\bold i(\bar X')}^*)$ are two constants and the solution $\bold{\bar x'}$ is the largest value of the function $f_{\bold i(\bar X')}$. Therefore, the solution $\bold{\bar x'}$ is also an optimal solution of
\begin{align*}	
	\max_{\bold{x} \in \mathcal{X}\cap \mathbb B^{n}}  \min_{i \in [m]} \frac{f_i(\bold{x})}{f_i(\bold{x}_i^*)}  	
	&=  \max \{ \eta : \eta \le \frac{f_{\bold i(\bar X')}(\bold{ x})}{f_{\bold i(\bar X')}(\bold{x}_{\bold i(\bar X')}^*)}, \bold{x} \in \mathcal{X}\cap \mathbb B^{n}, \eta \in \mathbb R \}.\\
	&=  \max \{ \eta : \eta \le \frac{f_{\bold i(\bar X')}(\bold{\bar  x'})}{f_{\bold i(\bar X')}(\bold{x}_{\bold i(\bar X')}^*)}, \eta \in \mathbb R \}.
\end{align*} 
\end{proof}

Proposition \ref{prop:no_needs_all_optimal} shows that solving  Problem \eqref{RSM} with certain $\alpha_i \neq f_{i}(\bold{x}_{i}^*)$ for all $i \in [m]$ may provide an optimal solution of  Problem \eqref{RSM3}. From Propositions \ref{Bound_RSM2} and \ref{prop:no_needs_all_optimal}, we arrive at a corollary for the final analysis of Problem \eqref{RSM3}.

\begin{corollary}
	From Propositions \ref{Bound_RSM2} and \ref{prop:no_needs_all_optimal}, for Problem \eqref{RSM} with $lb_i \le \alpha_i \le ub_i$ for all $i \in [m]$, an optimal solution $\bold{\bar x''}$ of Problem \eqref{RSM} is equivalent to an optimal solution of Problem \eqref{RSM3} if one of the following conditions holds. 
	\begin{itemize}		
	\item[(i)] if the solution $\bold{\bar x''}$ satisfies  $\min_{i \in [m]} \frac{f_i(\bold{\bar x''})}{ub_i} = \bar \eta_{relax}$, or
	\item[(ii)]  the condition of Proposition \ref{prop:no_needs_all_optimal} holds for $\bold{\bar x'} = \bold{\bar x''}$.
\end{itemize}	
\end{corollary}

{
	Finally, we derive a corollary that provides a strategy for the computational study.
	
	\begin{corollary}\label{prop:more_cuts}
		Given  $\bar X \subseteq V$, we have 
		$\max\{\eta:  (\eta,\bold{x}) \in \mathcal{C} \cap C (\bar X, \{\bold i(\bar X)\}), \ \bold{x} \in \mathcal{X}\cap \mathbb B^{n},  \eta \in \mathbb R\} \ge \max\{\eta: \ \forall i \in [m],\ (\eta,\bold{x}) \in \mathcal{C} \cap C (\bar X, \bold I(\bar X)), \ \bold{x} \in \mathcal{X}\cap \mathbb B^{n}, \ \eta \in \mathbb R\} \ge \max\{\eta: (\eta,\bold{x}) \in \mathcal{C} \cap C (\bar X, \bold I(\bar X)), \ \eta \le \frac{\theta_i}{\alpha_i}\ \ \forall i \in [m],\ (\theta_i, \bold{x}) \in  F_i, {i \in [m]}, \ \bold{x} \in \mathcal{X}\cap \mathbb B^{n}, \ \eta \in \mathbb R, \ \boldsymbol{\theta} \in \mathbb R^{m}\}$, where $\bold I(\bar X) = \{j \in [m]: \frac{f_j(\bar X)}{\alpha_j} = \frac{f_{\bold i(\bar X)}(\bar X)}{\alpha_{\bold i(\bar X)}}\}$, $\alpha_i = \max\{\theta_i: (\theta_i, \bold{x}) \in \bar F_i, \bold{x} \in \mathcal{X}\cap \mathbb B^{n}, \theta_i \in \mathbb R\}$ and $\bar F_i \supseteq F_i$ for all $i \in [m]$. 
			\end{corollary}
	\begin{proof}
		Given  $(\bar \eta,\boldsymbol{\bar  \theta},\bold{\bar x})$, where $\bar {\bold{x}} \in \mathbb B^{n} \cap \mathcal X \cap \mathcal C $, we have the following relations
		\begin{align*}
			\bar \eta & \le \min_{i \in \bold I(\bar X)}\frac{\bar \theta_i}{\alpha_i} \le \min_{i \in \bold I(\bar X)} \left\{ \frac{ f_{i}(\bar X) - \sum_{j\in \bar X\setminus \bar X}  \rho^{i}_j(V\setminus\{j\})+\sum_{j\in \bar X\setminus \bar X} \rho^{i}_j(\bar X) }{\alpha_{i}}  \right\}\\	
			& \le \frac{ f_{\bold i(\bar X)}(\bar X) - \sum_{j\in \bar X\setminus \bar X}  \rho^{\bold i(\bar X)}_j(V\setminus\{j\})+\sum_{j\in \bar X\setminus \bar X} \rho^{\bold i(\bar X)}_j(\bar X) }{\alpha_{\bold i(\bar X)}},
		\end{align*} 
		where the above relations follow from the definition of submodular inequality \eqref{sub_cut} and $\bold i(\bar X) \in \bold I(\bar X)$. This completes the proof.
	\end{proof}
	In Corollary \ref{coro:empty_set_cut}, we establish that the set of submodular inequalities $C (\emptyset, [m])$ satisfies the conditions of Corollary \ref{coro:sub_cuts_facet}. On the other hand, Corollary \ref{prop:more_cuts} shows that with the same RHS, adding a set of submodular inequalities  provides a tighter bound  compared to just adding one submodular inequality to RMP \eqref{RMP}, where the RHS is the value of $\frac{f_{\bold i(\bar X)}(\bar X)}{\alpha_{\bold i(\bar X)}} \le \frac{f_i(\bar X)}{\alpha_i}$ for a given $\bar X \in V$. Finally, Corollary \ref{prop:more_cuts} notes that as we solve a submodular maximization problem $\alpha_i = \max\{\theta_i: (\theta_i, \bold{x}) \in \bar F_i, \bold{x} \in \mathcal{X}\cap \mathbb B^{n}, \theta_i \in \mathbb R\}$, a subset of submodular inequalities defining the mixed-integer set $\bar F_i \supseteq F_i$ can be reused to derive a class of valid inequalities of RMP \eqref{RMP} for solving Problem \eqref{RSM3}. 
	}

In the next section, we design algorithms for solving Problem \eqref{RSM}, including the special case of Problem \eqref{RSM3}. The idea of Proposition \ref{Bound_RSM2} is that if we obtain all  the lower and upper bounds of the $m$ submodular maximization problems \eqref{SubMax}, we can calculate the optimality gap of the associated Problem \eqref{RSM3}.  Proposition \ref{Bound_RSM2} provides a strategy for solving Problem \eqref{RSM3}. That is, we could set  a time limit for each submodular maximization Problem \eqref{SubMax} to obtain upper and lower bounds of the problem. Using these bounds,  for all $i \in [m]$, we set  $\alpha_i$ as the lower bound $lb_i$. By solving the relaxation  $\max\{\eta:  (\eta, \bold{x}) \in \mathcal{C}, \ \bold{x} \in \mathcal{X}\cap \mathbb B^{n}, \ \eta \in \mathbb R$  we obtain the optimality gap of Problem \eqref{RSM3}.

\subsection{Algorithms}

In the final part of this section, we summarize the mentioned strategies and provide algorithms for Problem \eqref{RSM} including the special case of Problem \eqref{RSM3}, with $\boldsymbol{\alpha} = (f_1(\bold{x}_1^*),f_2(\bold{x}_2^*), \dots, f_m(\bold{x}_m^*))$. The core algorithm is a delayed constraint generation algorithm described in Algorithm \ref{alg:DCG}.  Algorithm \ref{alg:DCG} takes as input $\boldsymbol{\alpha}$, a subset of cuts defining $\mathcal C$ (could be empty), and a Boolean parameter $reduce$ that determines whether we consider  Proposition \ref{prop:MIP_valid} and Corollary \ref{prop:more_cuts}.  The True value of the parameter $reduce$ indicates that we consider the mixed-integer set $\mathcal F$ or $\mathcal F'$ that includes fewer submodular inequalities compared to adding all inequalities for $F_i, i \in  [m]$ under the False value of the parameter. The termination criteria can be a time limit, $T$ and/or an optimality gap tolerance, $\epsilon \in [0,1]$, where, for a lower bound on the optimal solution denoted as $\min(\Lambda)$ and an incumbent objective value $ \bar \eta$, the optimality gap is given by $ \bar \eta - \min(\Lambda)$. Note that the user can provide warm-start cuts for the set $\mathcal C$ of the RMP \eqref{RMP} as input. In particular, in Corollary \ref{coro:empty_set_cut}, we have shown that the set of submodular inequalities $C(\emptyset,[m])$ satisfies the facet conditions given in Corollary \ref{coro:sub_cuts_facet}. Therefore, in line \ref{algo:warmup} of Algorithm \ref{alg:DCG}, we add the facet-defining inequalities  $C(\emptyset,[m])$ to $\mathcal C$ as a class of warm-start cuts.

In line \ref{algo:getIncumbent} of the while loop of Algorithm \ref{alg:DCG}, we solve RMP \eqref{RMP} and get an incumbent solution ($\bar \eta, \bold{\bar x}$). In line \ref{algo:Omega}, based on the incumbent $\bold{\bar x}$, we form a set $\Lambda$ including $m$ values $\frac{f_i(\bold{\bar x})}{\alpha_i}$ for all $i \in [m]$. We compute $\min_{i \in [m]} \frac{f_i(\bold{\bar  x})}{\alpha_i}$ using the function $\min(\Lambda)$ that returns the minimal value of the elements of the set $\Lambda$. Note that the function $\min(\Lambda)$ provides a lower bound of Problem \eqref{RSM} based on the incumbent $\bold{\bar x}$. The lower bound is used to compute an optimality gap and obtain {the smallest value of the set $\left\{\frac{f_1(\bold{\bar x})}{\alpha_1}, \dots,\frac{f_i(\bold{\bar x})}{\alpha_i}, \dots, \frac{f_m(\bold{\bar x})}{\alpha_m} \right\}$, which is essential to determine the set $\mathcal F$ of Proposition \ref{prop:MIP_valid}.}  The for loop in  lines \ref{algo:add_cutsB} to \ref{algo:add_cutsE}  is for adding the submodular inequalities.  Given an incumbent $\bold{\bar x}$, if $reduce$ = True, the for loop adds fewer submodular inequalities to RMP \eqref{RMP} following Proposition \ref{prop:MIP_valid} and Corollary \ref{prop:more_cuts}, compared to the case that $reduce$ = False. 
Next, in  
Algorithm \ref{alg:FindSet}, we describe the separation routine, FindSetRoutine($\bold{\bar x}$,$i$), of this for loop. Recall that given an incumbent $\bar X$ and $i \in [m]$, Proposition \ref{prop:inter_cut} separates the incumbent $\bar X$ into two sets  $\tilde X_i \subseteq \bar X$ and $\bar X \setminus \tilde X_i \subseteq \bar X$, where $\tilde X_i \subseteq \bar X$ determines the set
$\mathcal{S}(i,\tilde X_i) = \{ j \in V\setminus \tilde X_i: \rho^{i}_j(\{k\}) = 0, \exists k \in \tilde X_i \}$. 
{Given an element $j \in V$, lines \ref{algoSep:for2B} to \ref{algoSep:for2E} first evaluate if there exists an element $k \in \bar X$ with $\rho^{i}_j(\{k\})=0$ for some $i \in [m]$. {Then, the algorithm determines if $j \in V$ can be a candidate of $S$, which is used to determine the set $\mathcal{S}(i,\tilde X_i)$, based on the condition shown in line \ref{algoSep:if3}, where $StopPt \in \mathbb N$ denotes the number of elements in $\bar X$ with zero marginal contribution.} }  If  $StopPt = 0$, then FindSetRoutine($\bold{\bar x}$,$i$) returns the original input $\bar X$. Here, line \ref{algoSep:if3} follows the condition \eqref{prop:inter_cut2} of Proposition \ref{prop:inter_cut} that allows us to consider the mixed-integer set $\mathcal F'$ as a valid set of submodular inequalities for  the set $\mathcal C$ of the RMP \eqref{RMP}.

\begin{algorithm}
	\SetAlgoLined
	 Input:  $\boldsymbol{\alpha} = (\alpha_1, \alpha_2, \dots, \alpha_m)$, $\mathcal C$, and a Boolean parameter $reduce$ \\
	{$\mathcal {C} \leftarrow \mathcal {C} \cap C(\emptyset,[m])$\\ \label{algo:warmup}}
{
	\While{Termination criteria not met}
	{
		Solve RMP \eqref{RMP} and obtain an incumbent ($\bar \eta, \bold{\bar x}$)\\\label{algo:getIncumbent}		
		\For{$i\in [m]$}
		{
			$\Lambda \leftarrow \Lambda \cup \left\{\frac{f_i(\bold{\bar x})}{\alpha_i}\right\}$\\\label{algo:Omega}
		}
			\For{$i\in [m]$}
			{\label{algo:add_cutsB}
				\If{$reduce =$ True}
				{ 
					\If{$\frac{f_i(\bold{\bar x})}{\alpha_i} = \min(\Lambda)$ and $\bar \eta > \frac{f_i(\bold{\bar x})}{\alpha_i}$}
					{ 
						$S \leftarrow$ FindSetRoutine($\bold{\bar x}$,$i$)\\
						Add a submodular inequality $\eta \le \frac{1}{ \alpha_i}(f_i(S) - \sum_{j\in S}  \rho^i_j(V\setminus\{j\})(1-x_j)+\sum_{j\in V\setminus S} \rho^i_j(S) x_j)$  to $\mathcal{C}$\\
					}
				}
				\Else 
				{
					\If{$\bar \eta > \frac{f_i(\bold{\bar x})}{\alpha_i}$}
					{ 
						$S \leftarrow$ FindSetRoutine($\bold{\bar x}$,$i$)\\
						Add a submodular inequality 				$\eta \le \frac{1}{ \alpha_i}(f_i(S) - \sum_{j\in S}  \rho^i_j(V\setminus\{j\})(1-x_j)+\sum_{j\in V\setminus S} \rho^i_j(S) x_j)$  to $\mathcal{C}$\\
					}			
				}
			}\label{algo:add_cutsE}
		
	}
}
	Return ($\bar \eta, \bold{\bar x}$) as the optimal value and solution.
	\caption{Delayed Constraint Generation Algorithm ($\boldsymbol{\alpha}$, $\mathcal C$, $reduce$)}
	\label{alg:DCG}
\end{algorithm}

\begin{algorithm}
	\SetAlgoLined
	Set  a stop point $StopPt \in \mathbb N$\\
	$Q \leftarrow \emptyset$ \\
	$S \leftarrow \emptyset$ \\
	\For{$j\in V$}
	{	
		\If{$\rho^{i}_j(\bar X)$ = 0}
		{
			$tmpQ \leftarrow Q$\\
			$counter \leftarrow 0$\\
			\For{$k\in \bar X$}
			{ \label{algoSep:for2B}
				\If{$\rho^{i}_j(\{k\})$ = 0}
				{
					$counter \leftarrow counter +1$\\
					$tmpQ \leftarrow tmpQ \cup \{k\}$\\
				}
				\If{$counter = StopPt$ and $f_i(tmpQ) = f_i(S \cup \{j\})+ \sum_{l \in tmpQ} \rho^{i}_l(S \cup\{j\})$}
				{\label{algoSep:if3}
					$S \leftarrow S \cup {j}$\\
					$Q \leftarrow Q \cup tmpQ$\\
				}
			}\label{algoSep:for2E}
		}
	}
	$S \leftarrow S \cup \{ \bar X \setminus Q\}$\\
	Return $S$
	\caption{FindSetRoutine($\bold{\bar x}$,$i$)}
	\label{alg:FindSet}
\end{algorithm}

 Finally, we present Algorithm \ref{alg:RSM2} for solving Problem \eqref{RSM3}, which is Problem \eqref{RSM} with a special choice of $\boldsymbol{\alpha} = (f_1(\bold{x}_1^*),f_2(\bold{x}_2^*), \dots, f_m(\bold{x}_m^*))$. Recall that  obtaining $\boldsymbol{\alpha} = (f_1(\bold{x}_1^*),f_2(\bold{x}_2^*), \dots, f_m(\bold{x}_m^*))$, before solving the corresponding RMP \eqref{RMP}, requires the solution of $m$ $\mathcal{NP}$-hard submodular maximization problems \eqref{SubMax}. However, even if we cannot solve the $m$ problems optimally, we can use Algorithm \ref{alg:RSM2} to find a feasible solution for RSM \eqref{RSM2} along with an optimality gap.  
Lines \ref{algRMP2:whileB} to \ref{algRMP2:whileE} follow a standard method for solving a submodular maximization problem using the submodular inequalities with some additional features. In lines \ref{algRMP2:inSubB} to \ref{algRMP2:inSubE}, when we finish solving a submodular maximization problem,  the lower and upper bounds are recorded. Furthermore, since the submodular inequalities for the corresponding submodular maximization problem can be reused for solving Problem \eqref{RSM}, we adapt and store the inequalities to the set $\mathcal {\bar C}$ for further usage in line \ref{algRMP2:inSubE}. After the for loop of Algorithm \ref{alg:RSM2}, we call Algorithm \ref{alg:DCG} based on a new vector $\boldsymbol{\alpha}$ and a set of warm-start cuts $\mathcal {\bar C}$. At the end of Algorithm \ref{alg:RSM2}, using the returned incumbent solution of Algorithm \ref{alg:DCG}, we are able to compute an optimality gap for Problem \eqref{RSM2}, where the computation of the gap follows Proposition  \ref{Bound_RSM2}.

\begin{algorithm}
	\SetAlgoLined
	Let $UB \leftarrow \infty$ be the upper bound of Problem \eqref{RSM3}\\
	Let $LB \leftarrow 0$ be the lower bound of Problem \eqref{RSM3}\\
	$SubCutReduction \leftarrow$ True\\
{
	\For{$i\in [m]$}
	{ \label{algRMP2:forB}
		Let $\bar F_i $ be a mixed-integer set derived from a subset of constraints for the $i$-th submodular maximization Problem \eqref{SubMax}  \\
		\While{$True$}
		{	\label{algRMP2:whileB}
			Solve a master problem $\max \{\eta: \bold{x} \in \mathcal{X} \cap  \mathbb B^n, \eta \in \mathbb R, (\eta,\bold{x})\in \bar F_i\}$ and get an incumbent ($\bar \eta, \bold{\bar x}$)	\\
			\If{Termination criteria met}
			{  \label{algRMP2:inSubB}
				$ub_i \leftarrow \bar \eta$\\
				$\bar {\alpha_i} \leftarrow lb_i \leftarrow f_i(\bold{\bar x})$\\
				Modify each submodular inequality of $\bar F_i$ to the form $\eta \le \frac{1}{\bar {\alpha_i}}(f_i(S) - \sum_{j\in S}  \rho^i_j(V\setminus\{j\})(1-x_j)+\sum_{j\in V\setminus S} \rho^i_j(S) x_j)$ 
				and add the modified inequalities to  $\mathcal {\bar C}$    	\\	\label{algRMP2:inSubE}		
				\bf{break}\;
							} 
			\Else
			{
				$S \leftarrow$ FindSetRoutine($\bold{\bar x}$,$i$)\\
				Add a submodular inequality $\eta \le f_i(S) - \sum_{j\in S}  \rho^i_j(V\setminus\{j\})(1-x_j)+\sum_{j\in V\setminus S} \rho^i_j(S) x_j $ to $\bar F_i$ \\
			}	
		} \label{algRMP2:whileE}				
	} \label{algRMP2:forE}

}
	$\boldsymbol{\bar \alpha} \leftarrow (\bar {\alpha_1}, \bar {\alpha_2}, \dots, \bar {\alpha_i})$\\
	($\bar \eta, \bold{\bar x}$) $\leftarrow$ Algorithm \ref{alg:DCG}($\boldsymbol{\bar \alpha}, \mathcal {\bar C},SubCutReduction$)\\ \label{algRMP2:callDCG}
	$UB \leftarrow \bar \eta$\\
	$LB \leftarrow \min_{i \in [m]} \frac{f_i(\bold{\bar x})}{ub_i}$\\
	$Gap \leftarrow \frac{UB-LB}{UB}$
	
	\caption{Solution Method for Problem \eqref{RSM3}}
	\label{alg:RSM2}
\end{algorithm}

\section{An Application on a Class of Water Sensor Placement Optimization Problems}\label{computation}

In this section, we apply the proposed algorithms to a class of sensor placement optimization problems in a water distribution network, where the goal of the deployed sensors is to detect contaminants in the network. Various objectives have been considered to quantify the effectiveness of the sensor deployment, such as the volume of the contaminated water  \citep{Kessler1998}, the contaminant   detection time   \citep{Kumar1997,Ostfeld2004}, or the population affected by the pollutants. We refer the reader to \cite{Berry2005,BWSN2008} for a detailed introduction to sensor placement optimization in real-world applications. In addition, \cite{MultiO,MultiO1,MultiO2,MultiO3,MultiO4,EN3} provide an introduction to multi-objective sensor placement optimization problems.

\subsection{A Model for Sensor Placement Optimization Problems}

In this subsection, we introduce the outbreak detection model of \cite{Leskovec07}  in a water distribution network \cite[see also,][]{KrauseWater2008}. 
Let $J$ be a set of possible contamination events corresponding to a source node $j\in J$ polluting the network with  probability  $p_j \in [0,1]$. Therefore, a network may have $|J|$ different contamination sources. Let $V$ be the set of all possible sensor locations and $S\subseteq V$ be a set of selected sensor placements. Note that each sensor $s\in S$ has its own cost $a_s \in \mathbb R_+$; the total cost $\sum_{s \in S}a_s$ of the selection $S$ must be  less than or equal to a given budget $b \in \mathbb R_+$. Let $W$ be a vector of edge flow velocities (time). 
Let $T(S,j)$ be a detection time that a set of sensors $S \subseteq V$ detects the contamination of a source $j \in J$. From the definition, the function  $T(\{s\},j)$ denotes the time that a sensor $s \in V$ detects the contamination of $j$. We then derive a relation   $T(S,j) = \min_{s \in S} T(\{s\},j)$ meaning the time for detecting a contamination of $j$ is the minimal time for the contamination detected by any sensor $s \in S$. Note that if the contamination of $j$ cannot be detected by the set $S$, the function $T(S,j)$ takes a  value of $\infty$. Following the definition of the detection time, we let $\beta_j(t)$ be a penalty function that denotes the amount of damage caused by a source $j \in J$ after a time $t$. 
Here, the amount of damage can be defined by users. For example, in a water distribution network, the associated damage can be the number of polluted nodes, the population affected by contamination, or the total cost of the contamination. For the case $t = \infty$, the function  $\beta_j(\infty)$ denotes the total amount of damage caused by the contamination of $j \in J$. Note that the penalty function is non-decreasing, where $\beta_j(t) \le \beta_j(t')$ for $t \le t'$ and $t,t' \in \mathbb R_+$.

Consider a water distribution network represented by a graph $G = (V,E,J)$, where $V$ is a set of nodes, $E$ is a set of directed edges, $J\subseteq V$ is a set of possible contamination sources. Based on the definition of the penalty function, given a set of sensors $S$ and a contamination source $j \in J$ in $G$, the penalty reduction is defined as $R_{G,W}(S,j) = \beta_j(\infty)-\beta_j(T(S,j))$. The penalty reduction measures  the amount of damage that can be avoided due to the contamination of $j$ after deploying a set $S$ of sensors
in the water distribution network. Recall that the probability of the event $j \in J$  is $p_j \in [0,1]$. For a set of possible contamination sources $J$ and a set of deployed sensors $S \subseteq V$, we consider the expected penalty reduction function $\mathcal{R}_{G,W}: 2^V \rightarrow \mathbb R$, where $\mathcal{R}_{G,W}(S) = \sum_{j\in J}p_jR_{G,W}(S,j)$ is  submodular  \cite[see,][for a proof of submodularity]{Leskovec07}. Below, we give an example to illustrate the outbreak detection model in a water distribution network.

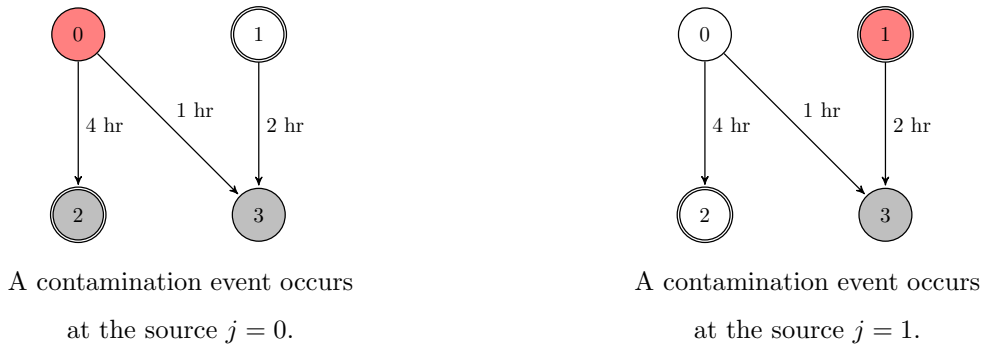
\begin{figure}[htbp]
	\begin{minipage}[t]{0.5\textwidth}
		\centering
		\scalebox{0.8}{\begin{tikzpicture}[->,>=stealth',shorten >=1pt,auto,node distance=3cm, semithick]
				\tikzstyle{every state}=[fill=none,draw=black,text=black]
				\node[state, accepting, fill=gray!50]  (2) []           {$2$};				
				\node[state, fill=red!50](0) [above of=2] {$0$};
				\node[state, accepting]  (1) [right of=0] {$1$};
				
				\node[state, fill=gray!50]  (3) [right of=2] {$3$};
				\path(0) edge   node {4 hr} (2);
				\path(0) edge   node {1 hr} (3);
				\path(1) edge   node {2 hr} (3);
				;
		\end{tikzpicture}}\\
		A contamination event occurs\\
		at the source $j=0$. 
	\end{minipage}
	\begin{minipage}[t]{0.5\textwidth}
		\centering	
		\scalebox{0.8}{\begin{tikzpicture}[->,>=stealth',shorten >=1pt,auto,node distance=3cm,
				semithick]
				\tikzstyle{every state}=[fill=none,draw=black,text=black]
				\node[state, accepting]                (2) []           {$2$};				
				\node[state](0) [above of=2] {$0$};
				\node[state, accepting, fill=red!50]   (1) [right of=0] {$1$};
				
				\node[state, fill=gray!50]  (3) [right of=2] {$3$};
				
				\path 
				(0) edge   node {4 hr} (2)
				(0) edge   node {1 hr} (3)
				(1) edge   node {2 hr} (3)
				;
		\end{tikzpicture}}\\
		A contamination event occurs\\ 
		at the source $j=1$.
	\end{minipage}
	\caption{An example introducing the penalty reduction via
		a network $G=(V,E,J)=(\{0,1,2,3\},\{(0,2),(0,3),(1,3)\},\{0,1\})$ with 4 nodes, 3 directed edges, and 2 possible contamination sources $J = \{0,1\}$, and $W=(4,1,2)$.}
	\label{fig:spex3}
\end{figure}

\begin{example}
Consider the water distribution network shown in Figure \ref{fig:spex3}. The network is represented as $G=(V,E,J)=(\{0,1,2,3\},\{(0,2),(0,3),(1,3)\},\{0,1\})$, where each directed edge $(i,j) \in E$ indicates  the water flow  from  node $i$ to  node $j$. Each edge $(i,j) \in E$ has a weight representing the flow time from $i$ to $j$. For example, the weight of the edge $(0,2)$ is 4, indicating that  it  takes 4 hours for  the water  to flow from  node $0$ to  node $2$. In the network $G$, we deploy a set of sensors $S=\{1,2\}$  indicated by two double circles on  nodes $1$ and $2$. For a set of two contamination sources $J=\{0,1\}$, we consider the following two cases in Figure \ref{fig:spex3}.

In the left subfigure of Figure \ref{fig:spex3}, a contamination event is indicated by the red node corresponding to  the contamination source $j=0$. Apart from the polluted red source node 0, the gray nodes 2 and 3, receiving the water flow from the source node $0$, are polluted if no sensor detects the contamination. Thus, the penalty function $\beta_0(\infty) = 3$ captures the number of polluted nodes without  any sensors, given by the two gray nodes and the red node. If a sensor is placed at  node 2 (indicated by the double circle), the contamination at  source node 0 will be detected after 4 hours; however, the sensor deployed at  node 1 cannot detect the contamination because there is no water flow from the source node $0$ to node $1$. Therefore, we conclude  $T(\{1,2\},0)=\min_{s \in \{1,2\}} T(\{s\},0)=4$, where $T(\{2\},0)=4$ and $T(\{1\},0) = \infty$. The associated penalty reduction $R_{G,W}(\{1,2\},0)=\beta_0(\infty)-\beta_0(T(\{1,2\},0))=3-2=1 $ denotes that under the contamination event at node $j = 0$, one node is not polluted because of the sensors deployed at $S$ in $G$. In other words, the set $S$ saves the damage to one node in $G$ under this contamination event.

In the right subfigure of Figure \ref{fig:spex3}, we consider another contamination source at $j = 1$. Two nodes (1 and 3) can be polluted by the water flow from  source $1$. However, since a sensor is placed at  source node 1, the  contamination event from node $1$ will be detected immediately. There are no nodes in $G$ polluted by the source $j = 1$ because of the sensor at  node 1. Thus, we conclude that the associated penalty reduction $R_{G,W}(\{1,2\},1)=\beta_1(\infty)-\beta_1(T(\{1,2\},1))=2-0=2 $, where $T(\{1,2\},0) = 0$.

Finally, for all $j \in J$, we assume that each contamination event of the source $j$ has the same probability $p_j=\frac{1}{2}$. The expected penalty reduction $\mathcal {R}_{G,W}(\{1,2\})=\sum_{j\in \{0,1\}}p_j R_{G,W}(\{1,2\},j)
=p_0R_{G,W}(\{1,2\},0)+p_1R_{G,W}(\{1,2\},1)
=0.5\times 1+0.5 \times 2=1.5 $.
\end{example}

Next, we formulate a robust variant of the outbreak detection problem, with uncertain water flow velocity along each edge. The uncertainty is due to hurricane disturbances, clogged pipes, and pump failures that may affect the flow velocity along the pipes. We  represent each scenario $i$  for  $i \in [m]$ with  $g_i = (G,W_i)$, where $W_i$ is a vector of velocities (weights) for the edges in $G$. Recall that $\mathcal{R}_{G,W}: 2^V \rightarrow \mathbb R$ is a submodular function. In our experiments, we let  $f_i = \mathcal{R}_{g_i}$ for all $i \in [m]$. Furthermore, the set of constraints $\mathcal X$ is given by $\{x:\sum_{i \in V}a_ix_i \le b\}$. Given a scenario $g_i$, the goal of the submodular maximization problem \eqref{SubMax} is to find an optimal solution that provides the maximal value of the expected penalty reduction for this scenario under the constraint $\sum_{i \in V}a_ix_i \le b$. In other words, for a scenario $g_i$ and the budget constraint, Problem \eqref{SubMax} aims to place a set of sensors that avoid the largest expected amount of damage (i.e.,  save the largest expected number of nodes) caused by contamination events in $J$. In contrast, given a set of sensors, Problem \eqref{RSM2} aims to find an optimal sensor placement that protects the largest expected amount of nodes in the worst case of $m$ scenarios. On the other hand, in Problem \eqref{RSM3}, given a scenario $g_i$ and a set of sensors $S$, we consider the proportion of the number of saved nodes by sensors  in $S$ to the maximal number of protected nodes with an optimal placement under scenario $i$, where the latter value is obtained by solving the $i$-th submodular maximization problem \eqref{SubMax}. 
In the following subsection, we  evaluate our proposed methods shown in Section \ref{section_methods} on real  water distribution networks.

\subsection{Computational Results} 
In this subsection, we report our computational experience with the proposed methods. We first introduce  the  three water distribution networks used in our computational study. We consider two networks, EN2 and EN3, from EPANET developed by the United States Environmental Protection Agency. Furthermore, we consider a network, BWSN1, from the battle of water sensor networks of \cite{BWSN2008}. Note that in a water distribution network, a facility, such as a junction, reservoir source, or tank, is represented by a node. A pipe is represented by a node pair $(i,j)$ denoting the direction of an edge from $i$ to $j$ (see Node1 and Node2 of PIPES in http://epanet.de/js/index.html.en). The network EN2 includes 36 nodes and 41 edges, EN3 includes 97 nodes and 117 edges, and BWSN1 includes 129 nodes and 168 edges.

Based on the three networks, we use the following parameters for Problem \eqref{RSM}. We set  the number of nodes to $|V| \in \{36,97,129\}$, where there are $|J| \in \{25,50\}$ contamination sources for EN3 and BWSN1, and $|J| \in \{12,25\}$ for the small-size network EN2. The probability of a contamination event  at a contamination source $j \in |J|$ is $p_j = \frac{1}{|J|}$. We generate $m \in \{50,100\}$ scenarios for each network, where the weights $W_i$ of the directed edges for a scenario $g_i$ are chosen  from a discrete uniform distribution $\mathcal U (1,10)$ for all $i \in [m]$. We consider a budget $b \in \{30,50\}$, where the cost of a sensor $a_i \in A$ is from a discrete uniform distribution $\mathcal U(5,10)$ for all $i \in V$. Note that given a fixed budget $b$, the different cost set $A$ may affect the number of sensors deployed in a network. For each setting $(|V|, b, m, |J|)$, we generate three instances  and report the average statistics. All algorithms are implemented in Python with Gurobi  8.1.1 Optimizer. We execute all experiments on a laptop with Intel Core i5-10210U 1.60 GHz CPU, 8 GB DRAM, x64 processor, and Windows 10 operating system. The time limit for each instance is set to 1800 seconds. We consider $\epsilon = 0$ and use the default integrality gap (MIPgap) of Gurobi, where a MIPgap of $10^{-4} \%$ is considered optimal.

First, we consider Problem \eqref{RSM2}, which is Problem \eqref{RSM} under the case $\boldsymbol{\alpha} = \bold{1}$. Algorithm \ref{alg:DCG} is used for solving the problem. The Baseline-RSM \eqref{RSM2} column provides baseline computational results for solving Problem \eqref{RSM2} with Algorithm  \ref{alg:DCG} using the parameters $reduce$= False, and with  $StopPt=0$ for the associated FindSetRoutine (Algorithm \ref{alg:FindSet}). That is,  Baseline-RSM \eqref{RSM2} with $reduce$= False considers all submodular inequalities for $F_i, i \in  [m]$ for $\mathcal C$ instead of considering  $\mathcal F$ with fewer submodular inequalities as shown in Proposition \ref{prop:MIP_valid}. Also, in  Baseline-RSM \eqref{RSM2}, since the parameter $StopPt$ of the associated FindSetRoutine is zero, given an incumbent solution $\bold{\bar x} \in \mathcal X$,  the algorithm does not utilize Proposition \ref{prop:inter_cut} that allows us to find a better set  than $\bar X$ to generate the corresponding submodular inequality. We consider two other  methods, Poly$\mathcal F$-Algo \ref{alg:DCG} and Poly$\mathcal F'$-Algo \ref{alg:DCG}, shown in the other two columns of Table \ref{Table:RSM_1} to evaluate the computational benefits of Propositions \ref{prop:MIP_valid} and \ref{prop:inter_cut} described in Section \ref{section_methods}, respectively. In Poly$\mathcal F$-Algo \ref{alg:DCG}, we consider Algorithm \ref{alg:DCG} with $reduce$ = True and the parameter $StopPt=0$ in the associated FindSetRoutine. That is,  RMP \eqref{RMP} uses the set $\mathcal F$ for deriving the cuts in $\mathcal C$.
In Poly$\mathcal F'$-Algo \ref{alg:DCG}, Problem \eqref{RMP} considers set $\mathcal F'$ shown in Proposition \ref{prop:inter_cut} for deriving the cuts in $\mathcal C$. Note that for set $\mathcal F'$, we let  $reduce$ = True and  $StopPt = 2$ in the FindSetRoutine of Algorithm \ref{alg:DCG}. 

We summarize our computational results in Table \ref{Table:RSM_1}. The Time-s column denotes the average computational time of three instances (in seconds). Note that the number in the parenthesis under the Time-s column denotes the number of instances that cannot be solved within the time limit of 1800 seconds. The average gap of the unsolved instances is reported in the Gap-\% column and we use a dash symbol to indicate when all three instances of each setting are solved optimally. The Iteration-\# column records the number of iterations to solve RMP \eqref{RMP}. The Cut-\# column reports the number of submodular inequalities added to  set $\mathcal C$ of RMP \eqref{RMP}. From the Time-s columns, we observe that Poly$\mathcal F$-Algo \ref{alg:DCG} is faster than the baseline, which demonstrates the effectiveness of Proposition \ref{prop:MIP_valid}. We note that the Baseline-RSM \eqref{RSM2} adds more inequalities to RMP \eqref{RMP}, leading to many unsolved instances for the EN3 and BWSN1 instances. Comparing Poly$\mathcal F$-Algo \ref{alg:DCG} and Poly$\mathcal F'$-Algo \ref{alg:DCG}, we observe that for most instances, Poly$\mathcal F'$-Algo \ref{alg:DCG} outperforms Poly$\mathcal F$-Algo \ref{alg:DCG} in both computational time and the number of added inequalities. This highlights the effectiveness of Proposition \ref{prop:inter_cut} in these instances.

\begin{sidewaystable}[htb]
	
	\caption{The computational results for RSM \eqref{RSM2}.}
	\label{Table:RSM_1}
	\begin{center}
		\scalebox{0.8}{\begin{tabular}{ 			 |p{1.5cm}p{0.6cm}p{0.6cm}p{0.6cm}p{0.6cm}
					||p{1.3cm}p{1.3cm}p{1.9cm}p{1.3cm}
					||p{1.3cm}p{1.3cm}p{1.9cm}p{1.3cm}
					||p{1.3cm}p{1.3cm}p{1.9cm}p{1.3cm}|}
				\hline				
				&  & & & & \multicolumn{4}{c||}{Baseline-RSM \eqref{RSM2}} &
				\multicolumn{4}{c||}{Poly$\mathcal F$-Algo \ref{alg:DCG}} &
				\multicolumn{4}{c|}{Poly$\mathcal F'$-Algo \ref{alg:DCG}} \\	
				
				\cline{6-17}				
				{Networks}&{$|V|$}&{$b$}&{$m$}&{$|J|$} 
				& {Time-s} & {Gap-\%} & {Iteration-\#} & {Cut-\#}
				& {Time-s} & {Gap-\%} & {Iteration-\#} & {Cut-\#}
				& {Time-s} & {Gap-\%} & {Iteration-\#} & {Cut-\#}\\

				\hline
\multirow{8}{6mm}{EN2} & \multirow{8}{6mm}{36}
& \multirow{2}{6mm}{50} 
& \multirow{4}{6mm}{100} & 25
& 1303(2) & 1.81 & 83 & 7439
& 691 & - & 171 & 635
& 114 & - & 75 & 244\\ 
~ & ~ & ~ & ~ & 12
& 286 & - & 32 & 3082
& 576 & - & 47 & 3551
& 69 & - & 14 & 1227\\ 
~ & ~ & \multirow{2}{6mm}{30} & ~ & 25 
& 971 & - & 57 & 4688
& 211 & - & 79 & 187
& 40 & - & 35 & 64\\ 
~ & ~ & ~ & ~ & 12
& 119 & - & 12 & 1111
& 158 & - & 16 & 850
& 47 & - & 11 & 626\\ 
~ & ~ & \multirow{2}{6mm}{50} & \multirow{4}{6mm}{50} & 25
& 693(2) & 2.45 & 110 & 4846
& 850 & - & 191 & 427
& 86 & - & 78 & 161\\ 
~ & ~ & ~ & ~ & 12
& 160 & - & 40 & 1913
& 173 & - & 33 & 1198
& 37 & - & 16 & 650\\ 
~ & ~ & \multirow{2}{6mm}{30} & ~ & 25 
& 510 & - & 56 & 2285
& 93 & - & 75 & 112
& 39 & - & 37 & 62\\ 
~ & ~ & ~ & ~ & 12
& 59 & - & 12 & 550
& 60 & - & 13 & 360
& 20 & - & 11 & 309\\ 
\hline

				\multirow{8}{6mm}{EN3} & \multirow{8}{6mm}{97}
				& \multirow{2}{6mm}{50} 
				& \multirow{4}{6mm}{100} & 50
				& (3) & 3.25 & 19 & 1880
				& (3) & 1.97 & 129 & 137
				& 1179(1) & 3.88 & 174 & 186\\ 
				~ & ~ & ~ & ~ & 25
				& (3) & 2.1 & 28 & 2557
				& 884(1) & 0.78 & 130 & 209
				& 355 & - & 106 & 176\\ 
				~ & ~ & \multirow{2}{6mm}{30} & ~ & 50
				& (3) & 5.74 & 27 & 2370
				& 610 & - & 106 & 140
				& 232 & - & 66 & 101 \\ 
				~ & ~ & ~ & ~ & 25
				& 756 & - & 24 & 1559
				& 117 & - & 41 & 76
				& 95 & - & 41 & 84 \\ 
				~ & ~ & \multirow{2}{6mm}{50} & \multirow{4}{6mm}{50} & 50
				& (3) & 3.72 & 32 & 1504
				& 1772(2) & 2.27 & 148 & 155
				& 767 & - & 159 & 165 \\ 
				~ & ~ & ~ & ~ & 25
				& (3) & 2.1 & 60 & 2325
				& 896 & - & 153 & 201
				& 280 & - & 117 & 156 \\ 
				~ & ~ & \multirow{2}{6mm}{30} & ~ & 50
				& 1422(1) & 1.5 & 49 & 1921
				& 329 & - & 93 & 94
				& 122 & - & 59 & 59 \\ 
				~ & ~ & ~ & ~ & 25 
				& 313 & - & 25 & 825
				& 84 & - & 42 & 76
				& 63 & - & 39 & 77 \\ 
				
				\hline			
\multirow{8}{6mm}{BWSN1} & \multirow{8}{6mm}{129} 
& \multirow{2}{6mm}{50} 
& \multirow{4}{6mm}{100} & 50
& (3) & 2.82 & 24 & 2397
& (3) & 2.39 & 162 & 180
& 1181(2) & 0.94 & 194 & 218\\ 
~ & ~ & ~ & ~ & 25 
& (3) & 2.27 & 46 & 3859
& 1465(1) & 3.12 & 189 & 258
& 173 & - & 66 & 83\\ 
~ & ~ & \multirow{2}{6mm}{30} & ~ & 50 
& 1155(1) & 1.4 & 25 & 1419
& 228 & - & 55 & 59
& 38 & - & 28 & 29\\ 
~ & ~ & ~ & ~ & 25
& 726 & - & 39 & 1500
& 206 & - & 61 & 69
& 38 & - & 35 & 38\\ 
~ & ~ & \multirow{2}{6mm}{50} & \multirow{4}{6mm}{50} & 50
& (3) & 3.61 & 36 & 1786
& (3) & 2.83 & 170 & 188
& 1572(1) & 0.1 & 223 & 246\\ 
~ & ~ & ~ & ~ & 25
& (3) & 4.18 & 57 & 2467
& 1261 & - & 180 & 220
& 93 & - & 62 & 76\\ 
~ & ~ & \multirow{2}{6mm}{30} & ~ & 50
& 1057 & - & 32 & 1098
& 133 & - & 45 & 48
& 21 & - & 23 & 25\\ 
~ & ~ & ~ & ~ & 25
& 316 & - & 28 & 707
& 39 & - & 39 & 45
& 17 & - & 29 & 31\\
				
\hline				
		\end{tabular}}
	\end{center}
\end{sidewaystable}

	\begin{sidewaystable}[htb]
	
	\caption{The computational results for RSM \eqref{RSM3}.}
	\label{Table:RSM3}
	\begin{center}
		\scalebox{0.8}{\begin{tabular}{ 			 |p{1.5cm}p{0.6cm}p{0.6cm}p{0.6cm}p{0.6cm}
					||p{1.3cm}p{1.3cm}p{1.9cm}p{1.3cm}
					||p{1.3cm}p{1.3cm}p{1.9cm}p{1.3cm}|}
				\hline				
				&  & & & & \multicolumn{4}{c||}{Baseline-RSM \eqref{RSM3} }  &
				\multicolumn{4}{c|}{Poly$\mathcal F'$-Algo \ref{alg:RSM2} with a finite $t$} \\	
				
				\cline{6-13}				
				{Networks}&{$|V|$}&{$b$}&{$m$}&{$|J|$} 
				& {Time-s} & {Gap-\%} & {Iteration-\#} & {Cut-\#}
				& {Time-s} & {Gap-\%} & {Iteration-\#} & {Cut-\#}\\
				\hline
				\multirow{8}{6mm}{EN2} & \multirow{8}{6mm}{36}
				& \multirow{2}{6mm}{50} 
				& \multirow{4}{6mm}{100} & 25
				& 1185 & - & 85 & 193 	
				& 984 & - & 2 & 3926 \\ 
				~ & ~ & ~ & ~ & 12
				& 581 & - & 31 & 137				
				& 386 & - & 1 & 2860 \\ %
				~ & ~ & \multirow{2}{6mm}{30} & ~ & 25
				& 595 & - & 33 & 131 				
				& 391 & - & 2 & 2067 \\ 
				~ & ~ & ~ & ~ & 12
				& 234 & - & 16 & 117
				& 190 & - & 1 & 1570 \\ 
				~ & ~ & \multirow{2}{6mm}{50} & \multirow{4}{6mm}{50} & 25 
				& 801 & - & 65 & 121
				& 583 & - & 2 & 2125 \\ %
				~ & ~ & ~ & ~ & 12
				& 272 & - & 23 & 80 
				& 215 & - & 1 & 1471 \\ 
				~ & ~ & \multirow{2}{6mm}{30} & ~ & 25
				& 340 & - & 32 & 80
				& 278 & - & 2 & 1047 \\ 
				~ & ~ & ~ & ~ &12
				& 148 & - & 16 & 68
				& 94 & - & 1 & 781 \\ 
				\hline

				\multirow{8}{6mm}{EN3} & \multirow{8}{6mm}{97}
				& \multirow{2}{6mm}{50} 
				& \multirow{4}{6mm}{100} & 50
				& (3) & N/A & N/A & N/A
				& 1613(2) & 3.42 & 12 & 2365 \\ 
				~ & ~ & ~ & ~ & 25
				& (3) & N/A & N/A & N/A
				& 1551 & 1.89 & 5 & 2843 \\ 
				~ & ~ & \multirow{2}{6mm}{30} & ~ & 50
				& (3) & N/A & N/A & N/A
				& 1142(1) & 1.13 & 3 & 2633 \\ 
				~ & ~ & ~ & ~ & 25
				& 1193 & - & 50 & 149 
				& 715 & - & 3 & 2329 \\ 
				~ & ~ & \multirow{2}{6mm}{50} & \multirow{4}{6mm}{50} & 50
				& (3) & N/A & N/A & N/A
				& 1632 & 1.45 & 9 & 1872 \\ 
				~ & ~ & ~ & ~ & 25
				& 1174(1) & N/A & 65 & 113
				& 970(1) & 1.21 & 6 & 2273 \\ 
				~ & ~ & \multirow{2}{6mm}{30} & ~ & 50
				& 1112 & - & 55 & 103
				& 588 & - & 5 & 1323 \\ 
				~ & ~ & ~ & ~ & 25
				& 605 & - & 38 & 86
				& 340 & - & 2 & 1137 \\ 
				\hline

				\multirow{8}{6mm}{BWSN1} & \multirow{8}{6mm}{129} 
				& \multirow{2}{6mm}{50} 
				& \multirow{4}{6mm}{100} & 50
				& (3) & N/A & N/A & N/A
				& 1638(2) & 1.9 & 15 & 3260 \\ 
				~ & ~ & ~ & ~ & 25
				& 1519(2) & N/A & 48 & 146
				& 1083 & - & 6 & 3245 \\ 
				~ & ~ & \multirow{2}{6mm}{30} & ~ & 50
				& 1477 & - & 25 & 123 
				& 961 & - & 3 & 1407 \\ 
				~ & ~ & ~ & ~ & 25
				& 1222 & - & 32 & 130 
				& 808 & - & 3 & 1760 \\ 
				~ & ~ & \multirow{2}{6mm}{50} & \multirow{4}{6mm}{50} & 50
				& (3) & N/A & N/A & N/A
				& 1438(2) & 1.03 & 13 & 2448 \\ 
				~ & ~ & ~ & ~ & 25
				& 1153 & - & 76 & 124
				& 1026 & - & 3 & 1639 \\ 
				~ & ~ & \multirow{2}{6mm}{30} & ~ & 50
				& 761 & - & 26 & 74
				& 416 & - & 1 & 719 \\ 
				~ & ~ & ~ & ~ & 25 
				& 555 & - & 29 & 77
				& 387 & - & 4 & 887 \\					
				
				\hline
				
		\end{tabular}}
	\end{center}
\end{sidewaystable}

Next, we consider RSM \eqref{RSM3}, which is RSM \eqref{RSM} under the case $\boldsymbol{\alpha} = (f_1(\bold{x}_1^*),f_2(\bold{x}_2^*), \dots, f_m(\bold{x}_m^*))$. From these previous experiments, we conclude that  using set $\mathcal F'$ in deriving the submodular inequalities is the best strategy for solving RSM \eqref{RSM}. Therefore, in Algorithm \ref{alg:RSM2}, we set  $reduce$ = True and  $StopPt = 2$ of the FindSetRoutine in Algorithm \ref{alg:DCG}. In these experiments, we aim to highlight the benefits of Algorithm \ref{alg:RSM2}. That is, we demonstrate different experiments on the $\bold{If}$-condition of lines \ref{algRMP2:inSubB}-\ref{algRMP2:inSubE} of Algorithm \ref{alg:RSM2}. {Baseline-RSM \eqref{RSM3} considers the basic method without any computational enhancements described in Section \ref{sec:AnaRSM3}. That is,}  algorithm \ref{alg:RSM2} without reusing the submodular inequalities generated for solving $m$ submodular maximization problems exactly to calculate $\alpha$. {For Baseline-RSM \eqref{RSM3}, we set  $\mathcal{\bar C} = \emptyset$ in line \ref{algRMP2:inSubE} and $t = \infty$ of Algorithm \ref{alg:RSM2}. Here,} the parameter $t = \infty$ indicates that Algorithm \ref{alg:RSM2} has to exactly compute $\boldsymbol{\alpha} = (f_1(\bold{x}_1^*),f_2(\bold{x}_2^*), \dots, f_m(\bold{x}_m^*))$ before solving RSM \eqref{RSM}.  Poly$\mathcal F'$-Algo \ref{alg:RSM2} with a finite $t$  demonstrates the effectiveness of Proposition \ref{Bound_RSM2}. That is, without completely solving $m$ submodular maximization problems, we aim to find a near-optimal solution with a provable optimality gap based on Proposition \ref{Bound_RSM2}. For Poly$\mathcal F'$-Algo \ref{alg:RSM2}, we set  $t = 15$s for $m = 100$ and $t = 30$s for $m = 50$. Note that because the time limit  is  1800 seconds, if $m$ submodular maximization problems take $t\times m$ seconds, then the time limit of algorithm \ref{alg:DCG} embedded in Algorithm \ref{alg:RSM2} is $1800-t\times m$ seconds.

Table \ref{Table:RSM3} provides the computational results of the three methods introduced in the previous paragraph. For {the instances that can be solved by both Baseline-RSM \eqref{RSM3} and Poly$\mathcal F'$-Algo \ref{alg:RSM2} with a finite $t$}, we observe that the setting $\mathcal{\bar C} = \emptyset$ slows down the performance of Algorithm \ref{alg:RSM2}. This shows the effectiveness of line \ref{algRMP2:inSubE} in Algorithm \ref{alg:RSM2}. We now consider the {the unsolved instances (N/A) of Table \ref{Table:RSM3}} and  observe that Poly$\mathcal F'$-Algo \ref{alg:RSM2} with a finite $t$ outperforms {Baseline-RSM \eqref{RSM3} significantly}. We note that in Baseline-RSM \eqref{RSM3}, there are many unsolved instances for EN3 and BWSN1, and the unsolved instances cannot provide a gap as indicated by the N/A symbol in Table \ref{Table:RSM3}. However, Poly$\mathcal F'$-Algo \ref{alg:RSM2} with a finite $t$ overcomes this issue and provides a small optimality gap for the instances unsolved within the time limit. Given that security  of the water distribution infrastructure is critical,   a high-quality  sensor deployment plan with a certifiable performance guarantee which is robust to disruptions as provided by Algorithm \ref{alg:RSM2} is highly desirable.

\section{Conclusion}\label{sec:conc}

We investigate mixed-integer programming methods and a polyhedral study for a class of robust submodular optimization problems. We start by introducing a fundamental robust submodular optimization problem, where the goal is to deal with
the worst case of a set of possible submodular functions. Several  propositions, including a facet condition on the submodular inequalities of the associated polyhedral structure, allow us to devise a delayed constraint generation method to solve the problem optimally. We also consider an extension of the fundamental robust submodular optimization problem that generalizes several robust submodular maximization subproblems of interest. For cases in which the submodular maximization subproblems cannot  be solved exactly within a  time limit, we provide a method for finding a feasible solution with a certifiable optimality gap. Our computational experiments on a sensor placement optimization problem for water distribution networks with real-world datasets demonstrate the effectiveness of the proposed methods.

\section*{Acknowledgments}
Simge K\"u\c{c}\"ukyavuz is supported, in part by, ONR Grant N00014-22-1-2602. Hao-Hsiang Wu is supported, in part by, NSTC Taiwan 111-2221-E-A49-079 and  109-2222-E-009-005-MY2. Hsin-Yi Huang is supported, in part by, NSTC Taiwan 109-2222-E-009-005-MY2.

	\bibliographystyle{apalike}
	\bibliography{general}

\begin{thebibliography}{}

\bibitem[Adibi et~al., 2022]{Adibi2022}
Adibi, A., Mokhtari, A., and Hassani, H. (2022).
\newblock Minimax optimization: The case of convex-submodular.
\newblock In {\em Proceedings of The 25th International Conference on
  Artificial Intelligence and Statistics}, volume 151, pages 3556--3580. PMLR.

\bibitem[Ahmed and Atamt\"urk, 2011]{Ahmed2011}
Ahmed, S. and Atamt\"urk, A. (2011).
\newblock Maximizing a class of submodular utility functions.
\newblock {\em Mathematical Programming}, 128(1):149--169.

\bibitem[Atamt\"urk, 2006]{Alper2006}
Atamt\"urk, A. (2006).
\newblock Strong formulations of robust mixed 0-1 programming.
\newblock {\em Mathematical Programming}, 108:235--250.

\bibitem[Atamt\"urk and G\'omez, 2020]{Alper2020}
Atamt\"urk, A. and G\'omez, A. (2020).
\newblock Submodularity in conic quadratic mixed 0-€"1 optimization.
\newblock {\em Operations Research}, 68(2):609--630.

\bibitem[Atamt\"urk and G\'omez, 2022]{Alper2022}
Atamt\"urk, A. and G\'omez, A. (2022).
\newblock Supermodularity and valid inequalities for quadratic optimization
  with indicators.
\newblock {\em Mathematical Programming}, pages 1--44.

\bibitem[Austin et~al., 2009]{EN3}
Austin, R., Choi, C., Peris, A., Ostfeld, A., and Lansey, K. (2009).
\newblock Multi-objective sensor placements with improved water quality models
  in a network with multiple junctions.
\newblock {\em World Environmental and Water Resources Congress 2009: Great
  Rivers}, pages 1--9.

\bibitem[Ben-Tal et~al., 2005]{BenTalSupply2005}
Ben-Tal, A., Golany, B., Nemirovski, A., and Vial, J.-P. (2005).
\newblock Retailer-supplier flexible commitments contracts: A robust
  optimization approach.
\newblock {\em Manufacturing and Service Operations Management}, 7(3):248--271.

\bibitem[Ben-Tal and Nemirovski, 1998]{Ben1998}
Ben-Tal, A. and Nemirovski, A. (1998).
\newblock Robust convex optimization.
\newblock {\em Mathematics of Operations Research}, 23(4):769--805.

\bibitem[Ben-Tal and Nemirovski, 1999]{Ben1999}
Ben-Tal, A. and Nemirovski, A. (1999).
\newblock Robust solutions of uncertain linear programs.
\newblock {\em Operations Research Letters}, 25:1--13.

\bibitem[Ben-Tal and Nemirovski, 2000]{Ben2000}
Ben-Tal, A. and Nemirovski, A. (2000).
\newblock Robust solutions of linear programming problems contaminated with
  uncertain data.
\newblock {\em Mathematical Programming}, 88:411--424.

\bibitem[Berry et~al., 2005]{Berry2005}
Berry, J., Fleischer, L., Hart, W.~E., Phillips, C.~A., and Watson, J. (2005).
\newblock Sensor placement in municipal water networks.
\newblock {\em Water Resources Planning and Management}, 131(3):237--243.

\bibitem[Bertsimas et~al., 2011]{Bertsimas2011}
Bertsimas, D., Brown, D.~B., and Caramanis, C. (2011).
\newblock Theory and applications of robust optimization.
\newblock {\em SIAM REVIEW}, 53(3):464--501.

\bibitem[Bertsimas et~al., 2013]{Bertsimas2013}
Bertsimas, D., Litvinov, E., Sun, X.~A., Zhao, J., and Zheng, T. (2013).
\newblock Adaptive robust optimization for the security constrained unit
  commitment problem.
\newblock {\em IEEE Transactions on Power Systems}, 28(1):52--63.

\bibitem[Bertsimas and Sim, 2003]{Bertsimas2003}
Bertsimas, D. and Sim, M. (2003).
\newblock Robust discrete optimization and network flows.
\newblock {\em Mathematical Programming}, 98:49--71.

\bibitem[Bertsimas and Sim, 2004]{Bertsimas2004}
Bertsimas, D. and Sim, M. (2004).
\newblock The price of robustness.
\newblock {\em Operations Research}, 52(1):35--53.

\bibitem[Bertsimas and Thiele, 2006]{BertsimasSupply2006}
Bertsimas, D. and Thiele, A. (2006).
\newblock A robust optimization approach to inventory theory.
\newblock {\em Operations Research}, 54(1):150--168.

\bibitem[Bogunovic et~al., 2017]{Bogunovic2017}
Bogunovic, I., Mitrović, S., Scarlett, J., and Cevher, V. (2017).
\newblock Robust submodular maximization: A non-uniform partitioning approach.
\newblock In {\em Proceedings of the 34th International Conference on Machine
  Learning}, ICML-2017, pages 508--516.

\bibitem[Boykov and Jolly, 2001]{Boykov2001}
Boykov, Y. and Jolly, M.-P. (2001).
\newblock Interactive graph cuts for optimal boundary and region segmentation
  of objects in n-d images.
\newblock In {\em Proceedings Eighth IEEE International Conference on Computer
  Vision}. IEEE.

\bibitem[Church and Velle, 1974]{MAXCOV74}
Church, R. and Velle, C.~R. (1974).
\newblock The maximal covering location problem.
\newblock {\em Papers in Regional Science}, 32(1):101--118.

\bibitem[Coniglio et~al., 2022]{Coniglio2022}
Coniglio, S., Furini, F., and Ljubi\'c‡, I. (2022).
\newblock Submodular maximization of concave utility functions composed with a
  set-union operator with applications to maximal covering location problems.
\newblock {\em Mathematical Programming}, 196:9--56.

\bibitem[Cordeau et~al., 2019]{Cordeau2019}
Cordeau, J.-F., Furini, F., and Ljubić, I. (2019).
\newblock Benders decomposition for very large scale partial set covering and
  maximal covering location problems.
\newblock {\em European Journal of Operational Research}, 275(3):882--896.

\bibitem[Dorini et~al., 2004]{Ostfeld2004}
Dorini, G., Jonkergouw, P., Kapelan, Z., and di~Pierro, F. (2004).
\newblock An efficient algorithm for sensor placement in water distribution
  systems.
\newblock {\em Journal of Water Resources Planning and Management},
  130(5):377--385.

\bibitem[Dorini et~al., 2006]{MultiO4}
Dorini, G., Jonkergouw, P., Kapelan, Z., and di~Pierro, F. (2006).
\newblock An efficient algorithm for sensor placement in water distribution
  systems.
\newblock In {\em Eighth Annual Water Distribution Systems Analysis Symposium},
  WDSA.

\bibitem[Feige, 1998]{Feige1998}
Feige, U. (1998).
\newblock A threshold of ln n for approximating set cover.
\newblock {\em Journal of the ACM}, 45(4):634--652.

\bibitem[Feige et~al., 2011]{Feige2011}
Feige, U., Mirrokni, V.~S., and Vondr\'ak, J. (2011).
\newblock Maximizing non-monotone submodular functions.
\newblock {\em SIAM Journal on Computing}, 40(4):1133--1153.

\bibitem[Fischetti et~al., 2018]{FischettiSocial2018}
Fischetti, M., Kahr, M., Leitner, M., Monaci, M., and Ruthmair, M. (2018).
\newblock Least cost influence propagation in (social) networks.
\newblock {\em Mathematical Programming}, 170:293--325.

\bibitem[Ghaoui et~al., 2003]{Ghaoui2003}
Ghaoui, L.~E., Oks, M., and Oustry, F. (2003).
\newblock Worst-case value-at-risk and robust portfolio optimization: A conic
  programming approach.
\newblock {\em Operations Research}, 51(4):543--556.

\bibitem[Goldfarb and Iyengar, 2003]{Goldfarb2003}
Goldfarb, D. and Iyengar, G. (2003).
\newblock Robust portfolio selection problems.
\newblock {\em Mathematics of Operations Research}, 28(1):1--38.

\bibitem[Golovin and Krause, 2011]{Golovin2011}
Golovin, D. and Krause, A. (2011).
\newblock Adaptive submodularity: Theory and applications in active learning
  and stochastic optimization.
\newblock {\em Journal of Artificial Intelligence Research}, 42:427--486.

\bibitem[G\'omez, 2018]{Andres2018}
G\'omez, A. (2018).
\newblock Submodularity and valid inequalities in nonlinear optimization with
  indicator variables.
\newblock {\em Optimization-Online}.

\bibitem[G\"unne\c{c} et~al., 2019]{Dilek2019}
G\"unne\c{c}, D., Raghavan, S., and Zhang, R. (2019).
\newblock Least-cost influence maximization on social networks.
\newblock {\em INFORMS Journal on Computing}, 32(2):289--302.

\bibitem[He and Kempe, 2016]{RIM2016}
He, X. and Kempe, D. (2016).
\newblock Robust influence maximization.
\newblock In {\em Proceedings of the 22nd ACM SIGKDD International Conference
  on Knowledge Discovery and Data Mining}, KDD '16, pages 885--894, New York,
  NY, USA. ACM.

\bibitem[Huang et~al., 2006]{MultiO1}
Huang, J.~J., McBean, E.~A., and James, W. (2006).
\newblock Multiobjective optimization for monitoring sensor placement in water
  distribution systems.
\newblock In {\em Eighth Annual Water Distribution Systems Analysis Symposium},
  WDSA.

\bibitem[Jegelka and Bilmes, 2011]{Jegelka2011}
Jegelka, S. and Bilmes, J. (2011).
\newblock Submodularity beyond submodular energies: Coupling edges in graph
  cuts.
\newblock In {\em CVPR}. IEEE.

\bibitem[Jiang et~al., 2012]{Jiang2012}
Jiang, R., Zhang, M., Li, G., and Guan, Y. (2012).
\newblock Benders decomposition for the two-stage security constrained robust
  unit commitment problem.
\newblock In {\em Proceedings of IIE Annual Conference}, pages 1--10.

\bibitem[Kempe et~al., 2003]{KKT03}
Kempe, D., Kleinberg, J., and Tardos, {\'E}. (2003).
\newblock Maximizing the spread of influence through a social network.
\newblock In {\em Proceedings of the Ninth ACM SIGKDD International Conference
  on Knowledge Discovery and Data Mining}, KDD '03, pages 137--146, New York,
  NY, USA. ACM.

\bibitem[Kessler et~al., 1998]{Kessler1998}
Kessler, A., Ostfeld, A., and Sinai, G. (1998).
\newblock Detecting accidental contaminations in municipal water networks.
\newblock {\em Journal of Water Resources Planning and Management},
  124(4):192--198.

\bibitem[K{\i}l{\i}n{\c{c}}-Karzan et~al., 2022]{Fatma2022}
K{\i}l{\i}n{\c{c}}-Karzan, F., K\"u\c{c}\"ukyavuz, S., and Lee, D. (2022).
\newblock Joint chance-constrained programs and the intersection of mixing sets
  through a submodularity lens.
\newblock {\em Mathematical Programming}, 195:283--326.

\bibitem[K{\i}l{\i}n{\c{c}}-Karzan et~al., 2020]{kilincc2020conic}
K{\i}l{\i}n{\c{c}}-Karzan, F., K{\"u}{\c{c}}{\"u}kyavuz, S., Lee, D., and
  Shafieezadeh-Abadeh, S. (2020).
\newblock Conic mixed-binary sets: Convex hull characterizations and
  applications.
\newblock {\em ArXiv:2012.14698}.

\bibitem[Kouvelis and Yu, 1997]{Kouvelis1997}
Kouvelis, P. and Yu, G. (1997).
\newblock Robust discrete optimization and its applications.
\newblock {\em Kluwer Academic Publishers, Norwell, MA,}.

\bibitem[Krause and Golovin, 2012]{SubMaxApp2012}
Krause, A. and Golovin, D. (2012).
\newblock Submodular function maximization.
\newblock {\em Tractability: Practical Approaches to Hard Problems}, 3(19).

\bibitem[Krause et~al., 2008a]{KrauseWater2008}
Krause, A., Leskovec, J., Guestrin, C., VanBriesen, J., and Faloutsos, C.
  (2008a).
\newblock Efficient sensor placement optimization for securing large water
  distribution networks.
\newblock {\em Journal of Water Resources Planning and Management},
  134(6):516--526.

\bibitem[Krause et~al., 2008b]{RSM2008}
Krause, A., McMahan, H.~B., Guestrin, C., and Gupta, A. (2008b).
\newblock Robust submodular observation selection.
\newblock {\em Journal of Machine Learning Research}, 9(93):2761--2801.

\bibitem[Krause et~al., 2008c]{KrauseSensor2008}
Krause, A., Singh, A., and Guestrin, C. (2008c).
\newblock Near-optimal sensor placements in {Gaussian} processes: Theory,
  efficient algorithms and empirical studies.
\newblock {\em Journal of Machine Learning Research}, 9:235--284.

\bibitem[K\"u\c{c}\"ukyavuz and Yu, 2023]{KYTut2023}
K\"u\c{c}\"ukyavuz, S. and Yu, Q. (2023).
\newblock Mixed-integer programming approaches to generalized submodular
  optimization and its applications.
\newblock {\em INFORMS Tutorials in Operations Research}.
\newblock ArXiv:2304.00479.

\bibitem[Kumar et~al., 1997]{Kumar1997}
Kumar, A., Kansal, M.~L., and Arora, G. (1997).
\newblock Identification of monitoring stations in water distribution system.
\newblock {\em Journal of Environmental Engineering}, 123(8):746--752.

\bibitem[Leskovec et~al., 2007]{Leskovec07}
Leskovec, J., Krause, A., Guestrin, C., Faloutsos, C., VanBriesen, J., and
  Glance, N. (2007).
\newblock Cost-effective outbreak detection in networks.
\newblock In {\em Proceedings of the 13th ACM SIGKDD international conference
  on Knowledge discovery and data mining}, KDD '07, pages 420--429, New York,
  NY, USA. ACM.

\bibitem[Lin and Bilmes, 2011]{Lin2011}
Lin, H. and Bilmes, J. (2011).
\newblock A class of submodular functions for document summarization.
\newblock In {\em Proceedings of the 49th Annual Meeting of the Association for
  Computational Linguistics: Human Language Technologies}, pages 510--520.
  Association for Computational Linguistics.

\bibitem[Mulvey et~al., 1995]{Mulvey1995}
Mulvey, J.~M., Vanderbei, R.~J., and Zenios, S.~A. (1995).
\newblock Robust optimization of large-scale systems.
\newblock {\em Operations Research}, 43(2):264--281.

\bibitem[Nannicini et~al., 2019]{Nannicini2019}
Nannicini, G., Sartor, G., Traversi, E., and Wolfler-Calvo, R. (2019).
\newblock An exact algorithm for robust influence maximization.
\newblock In {\em Integer Programming and Combinatorial Optimization}, page
  313–326. Lecture Notes in Computer Science.

\bibitem[Nemhauser and Wolsey, 1981]{NW81}
Nemhauser, G. and Wolsey, L. (1981).
\newblock Maximizing submodular set functions: Formulations and analysis of
  algorithms.
\newblock In Hansen, P., editor, {\em Annals of Discrete Mathematics (11)
  Studies on Graphs and Discrete Programming}, volume~59 of {\em North-Holland
  Mathematics Studies}, pages 279 -- 301. North-Holland Mathematics Studies.

\bibitem[Nemhauser et~al., 1978]{NWF78}
Nemhauser, G., Wolsey, L., and Fisher, M. (1978).
\newblock An analysis of approximations for maximizing submodular set
  functions---{I}.
\newblock {\em Mathematical Programming}, 14(1):265--294.

\bibitem[Orlin et~al., 2018]{Orlin2018}
Orlin, J.~B., Leskovec, J., and Udwani, R. (2018).
\newblock Robust monotone submodular function maximization.
\newblock {\em Mathematical Programming}, 172:505--537.

\bibitem[{Ostfeld et al.}, 2008]{BWSN2008}
{Ostfeld et al.} (2008).
\newblock The battle of the water sensor networks (bwsn): A design challenge
  for engineers and algorithms.
\newblock {\em Journal of Water Resources Planning and Management},
  134(6):556--568.

\bibitem[Powers et~al., 2016]{Powers2016}
Powers, T., Bilmes, J., Wisdom, S., Krout, D.~W., and Atlas, L. (2016).
\newblock Constrained robust submodular optimization.
\newblock In {\em Advances in Neural Information Processing Systems 30},
  NeurIPS-2016.

\bibitem[Preis and Ostfeld, 2006]{MultiO2}
Preis, A. and Ostfeld, A. (2006).
\newblock Multiobjective sensor design for water distribution systems security.
\newblock In {\em Eighth Annual Water Distribution Systems Analysis Symposium},
  WDSA.

\bibitem[Shen and Jiang, 2023]{Shen2022}
Shen, H. and Jiang, R. (2023).
\newblock Chance-constrained set covering with {Wasserstein} ambiguity.
\newblock {\em Mathematical Programming}, 198:621--674.

\bibitem[Shi et~al., 2022]{Shi2022}
Shi, X., Prokopyev, O.~A., and Zeng, B. (2022).
\newblock Sequence independent lifting for a set of submodular maximization
  problems.
\newblock {\em Mathematical Programming}, 196:69--€"114.

\bibitem[Staib et~al., 2019]{Staib2019}
Staib, M., Wilder, B., and Jegelka, S. (2019).
\newblock Distributionally robust submodular maximization.
\newblock In {\em Proceedings of the 22nd International Conference on
  Artificial Intelligence and Statistics}, volume~89, pages 506--516. PMLR.

\bibitem[T\"ut\"unc\"u and Koenig, 2004]{Tutuncu2004}
T\"ut\"unc\"u, R. and Koenig, M. (2004).
\newblock Robust asset allocation.
\newblock {\em Annals of Operations Research}, 132:157--187.

\bibitem[Watson et~al., 2004]{MultiO}
Watson, J.-P., Greenberg, H.~J., and Hart, W.~E. (2004).
\newblock A multiple-objective analysis of sensor placement optimization in
  water networks.
\newblock In {\em World Water and Environmental Resources Congress}.

\bibitem[Wu and K\"{u}\c{c}\"{u}kyavuz, 2018]{first2016}
Wu, H. and K\"{u}\c{c}\"{u}kyavuz, S. (2018).
\newblock A two-stage stochastic programming approach for influence
  maximization in social networks.
\newblock {\em Computational Optimization and Applications}, 69(3):563--595.

\bibitem[Wu and K\"{u}\c{c}\"{u}kyavuz, 2019]{Wu2017}
Wu, H. and K\"{u}\c{c}\"{u}kyavuz, S. (2019).
\newblock Probabilistic partial set covering with an oracle for chance
  constraints.
\newblock {\em SIAM Journal on Optimization}, 29(1):690--718.

\bibitem[Wu and Walski, 2006]{MultiO3}
Wu, Z.~Y. and Walski, T. (2006).
\newblock Multiobjective optimization of sensor placement in water distribution
  systems.
\newblock In {\em Eighth Annual Water Distribution Systems Analysis Symposium},
  WDSA.

\bibitem[Yu and Ahmed, 2017]{Yu2017}
Yu, J. and Ahmed, S. (2017).
\newblock Maximizing a class of submodular utility functions with constraints.
\newblock {\em Mathematical Programming}, 162(1-2):145--164.

\bibitem[Yu and K\"u\c{c}\"ukyavuz, 2021a]{QYuk2021}
Yu, Q. and K\"u\c{c}\"ukyavuz, S. (2021a).
\newblock An exact cutting plane method for $k$-submodular function
  maximization.
\newblock {\em Discrete Optimization}, 42:100670.

\bibitem[Yu and K\"u\c{c}\"ukyavuz, 2021b]{QYu2021}
Yu, Q. and K\"u\c{c}\"ukyavuz, S. (2021b).
\newblock A polyhedral approach to bisubmodular function minimization.
\newblock {\em Operations Research Letters}, 49(1):5--10.

\bibitem[Yu and K\"u\c{c}\"ukyavuz, 2022]{YuQ2022}
Yu, Q. and K\"u\c{c}\"ukyavuz, S. (2022).
\newblock On constrained mixed-integer {DR}-submodular minimization.
\newblock {\em Arxiv:2211.07726}.

\bibitem[Yu and K\"u\c{c}\"ukyavuz, 2023]{QYu2023}
Yu, Q. and K\"u\c{c}\"ukyavuz, S. (2023).
\newblock Strong valid inequalities for a class of concave submodular
  minimization problems under cardinality constraints.
\newblock {\em Mathematical Programming}, pages 1--59.

\bibitem[Zeng and Zhao, 2013]{Zeng2013}
Zeng, B. and Zhao, L. (2013).
\newblock Solving two-stage robust optimization problems using a
  column-and-constraint generation method.
\newblock {\em Operations Research Letters}, 41:457--461.

\bibitem[Zhao and Zeng, 2012]{Zhao2012}
Zhao, L. and Zeng, B. (2012).
\newblock Robust unit commitment problem with demand response and wind energy.
\newblock In {\em Proceedings of IEEE Power and Energy Society General
  Meeting}, pages 1--8. IEEE.

\bibitem[Zheng et~al., 2019]{ZhengIISE2019}
Zheng, K., Albert, L.~A., Luedtke, J.~R., and Towle, E. (2019).
\newblock A budgeted maximum multiple coverage model for cybersecurity planning
  and management.
\newblock {\em IISE Transactions}, 51(12):1303--1317.

\end{thebibliography}

\end{document}